\documentclass[10pt,reqno]{amsart}
\usepackage{amssymb}
\usepackage{amsfonts}
\usepackage{amsthm}
\usepackage{amsmath}

\newcommand{\Reals}{{\mathbb{R}}}
\newcommand{\Cmplx}{{\mathbb{C}}}
\newcommand{\Ints}{{\mathbb{Z}}}

\newcommand{\IC}{{\mathbb{C}}}

\newcommand{\ID}{{\mathbb{D}}}

\newcommand{\Disk}{{\mathbb{D}}}

\newcommand{\RR}{{{}_{R}}}

\newcommand{\g}{{\mathfrak{g}}}
\newcommand{\A}{{\mathfrak{a}}}
\newcommand{\U}{{\mathbb{U}}}
\newcommand{\LL}{{\mathfrak{l}}}
\newcommand{\LLL}{{\mathbb{L}}}

\newcommand{\C}{{\mathcal{C}}}

\let\det=\undefined\DeclareMathOperator*{\det}{det}
\let\Re=\undefined\DeclareMathOperator*{\Re}{Re}
\let\Im=\undefined\DeclareMathOperator*{\Im}{Im}

\DeclareMathOperator*{\tr}{tr} \DeclareMathOperator*{\diag}{diag}
\DeclareMathOperator*{\supp}{supp} \DeclareMathOperator*{\Id}{Id}

\newtheorem{theorem}{Theorem}[section]
\newtheorem{prop}[theorem]{Proposition}
\newtheorem{lemma}[theorem]{Lemma}
\newtheorem{coro}[theorem]{Corollary}
\theoremstyle{definition}
\newtheorem{definition}[theorem]{Definition}
\newtheorem{remark}[theorem]{Remark}

\newcounter{smalllist}
\newenvironment{SL}{\begin{list}{{\rm(\roman{smalllist})}}{%
\setlength{\topsep}{0mm}\setlength{\parsep}{0mm}\setlength{\itemsep}{0mm}%
\setlength{\labelwidth}{3em}\setlength{\leftmargin}{3em}\usecounter{smalllist}%
}}{\end{list}}

\begin{document}

\title{CMV: the unitary analogue of Jacobi matrices}
\author[R.~Killip and I.~Nenciu]{Rowan Killip$^1$ and Irina Nenciu}
\address{Rowan Killip\\
         UCLA Mathematics Department\\
         Box 951555\\
         Los Angeles, CA 90095}
\email{killip@math.ucla.edu}
\thanks{$^1$ Supported in part by NSF grant DMS-0401277 and a Sloan Foundation Fellowship.}

\address{Irina Nenciu\\
         Mathematics 253-37\\
         Caltech\\
         Pasadena, CA 91125}
\email{nenciu@caltech.edu}

\date{\today}

\dedicatory{Happy 60th birthday, Percy Deift.}

\begin{abstract}

We discuss a number of properties of CMV matrices, by which we mean
the class of unitary matrices recently introduced by Cantero, Moral,
and Velazquez.  We argue that they play an equivalent role among
unitary matrices to that of Jacobi matrices among all Hermitian
matrices. In particular, we describe the analogues of well-known
properties of Jacobi matrices: foliation by co-adjoint orbits, a
natural symplectic structure, algorithmic reduction to this shape,
Lax representation for an integrable lattice system
(Ablowitz-Ladik), and the relation to orthogonal polynomials.

As offshoots of our analysis, we will construct action/angle
variables for the finite Ablowitz-Ladik hierarchy and describe the
long-time behaviour of this system.
\end{abstract}

\maketitle

\section{Introduction}

For many reasons, it is natural to
regard Jacobi matrices as occupying a certain privileged place
among all Hermitian matrices.  By `Jacobi matrix' we mean a
tri-diagonal matrix
\begin{equation}\label{Jmat}
J=\begin{bmatrix}
 b_1  &  a_1 &        & \\
 a_1  &  b_2 & \ddots & \\
      &\ddots& \ddots & a_{n-1}\\
      &      & a_{n-1}&  b_n
\end{bmatrix}
\end{equation}
with $a_j>0$, $b_j\in\Reals$.  The main purpose of this paper is
explain why we consider a family of unitary matrices introduced
recently by Cantero, Moral, and Velazquez (see \cite{CMV}) as
playing the corresponding role among unitary matrices.

\begin{definition}\label{D:CMV}
Given coefficients $\alpha_0,\ldots,\alpha_{n-2}$ in $\Disk$ and
$\alpha_{n-1}\in S^1$, let $\rho_k=\sqrt{1-|\alpha_k|^2}$, and
define $2\times 2$ matrices
$$
\Xi_k = \begin{bmatrix} \bar\alpha_k & \rho_k \\ \rho_k & -\alpha_k
\end{bmatrix}
$$
for $0\leq k\leq n-2$, while $\Xi_{-1}=[1]$ and
$\Xi_{n-1}=[\bar\alpha_{n-1}]$ are $1 \times1$ matrices. From these,
form the $n\times n$  block-diagonal matrices
$$
\mathcal{L}=\diag\bigl(\Xi_0   ,\Xi_2,\Xi_4,\ldots\bigr)
\quad\text{and}\quad
\mathcal{M}=\diag\bigl(\Xi_{-1},\Xi_1,\Xi_3,\ldots\bigr).
$$
The \emph{CMV matrix} associated to the coefficients
$\alpha_0,\ldots,\alpha_{n-1}$ is $\C=\mathcal{LM}$.
\end{definition}

Following \cite{Simon1}, we will refer to the numbers $\alpha_k$ as Verblunsky coefficients .
A related system of matrices was discovered independently by Tao and Thiele, \cite{TaoThiele},
in connection with the non-linear Fourier transform.  Their matrices are bi-infinite and
correspond to setting odd-indexed Verblunsky coefficients to zero.  (This has the effect of
doubling the spectrum; cf. \cite[p. 84]{Simon1}.)

Expanding out the matrix product $\mathcal{L}\mathcal{M}$ is rather labourious.
In Figure~\ref{F1}, we show the result for $n=8$.  As can be seen from this example,
CMV matrices have a rather rigid structure:

\begin{definition}\label{D:CMVshape}
We say that an $n\times n$ matrix has \emph{CMV shape} if the
entries have the following pattern of horizontal $2\times 4$ blocks:
\begin{equation}\label{E:CMVshape}
\begin{bmatrix}
* & * & + & 0 & 0 & 0 & 0 & 0 \\
+ & * & * & 0 & 0 & 0 & 0 & 0 \\
0 & * & * & * & + & 0 & 0 & 0 \\
0 & + & * & * & * & 0 & 0 & 0 \\
0 & 0 & 0 & * & * & * & + & 0 \\
0 & 0 & 0 & + & * & * & * & 0 \\
0 & 0 & 0 & 0 & 0 & * & * & * \\
0 & 0 & 0 & 0 & 0 & + & * & *
\end{bmatrix}
\text{ or }
\begin{bmatrix}
* & * & + & 0 & 0 & 0 & 0 \\
+ & * & * & 0 & 0 & 0 & 0 \\
0 & * & * & * & + & 0 & 0 \\
0 & + & * & * & * & 0 & 0 \\
0 & 0 & 0 & * & * & * & + \\
0 & 0 & 0 & + & * & * & * \\
0 & 0 & 0 & 0 & 0 & * & *
\end{bmatrix}
\end{equation}
where $+$ represents a positive entry and $*$ represents a possibly
non-zero entry.  The top left corner of a CMV matrix always has the
$2\times 3$ structure depicted above, whereas the bottom right
corner will consist of a $2\times3$ block if $n$ is even (left
matrix above) or a $1\times 2$ block (right matrix above) if $n$ is
odd.
\end{definition}

Naturally, CMV matrices have CMV shape; a proof of this can be found in the original
paper of Cantero, Moral, and Vel\'azquez, \cite{CMV}.
Conversely, \cite{CMV2} shows that a unitary matrix with CMV shape must be a CMV matrix.
(Alternate proofs of these results can be found at the end of Section~\ref{s3}.)
This will be very useful for us since shape and unitarity
properties are easier to check than comparing all matrix entries.
To facilitate describing the shape of matrices we
will use the following terms:

\begin{definition}
The \emph{upper staircase} of a matrix consists of the first non-zero entry in every column (or last in
every row).  Conversely, the \emph{lower staircase} consists of the first non-zero entries in the rows.
\end{definition}

The entries marked $+$ are particularly important so we give them a name too.

\begin{definition}\label{D:Exposed}
The entries marked $+$ are precisely $(2,1)$ and those of the form
$(2j-1,2j+1)$ and $(2j+2,2j)$ with $j\geq 1$.  We will refer to
these as the \emph{exposed} entries of the CMV matrix.
\end{definition}

\begin{figure}
\begin{equation*}\label{F:CMV}
\begin{bmatrix}
\bar \alpha_0 & \rho_0\bar\alpha_1 & \rho_0\rho_1          & 0 & 0 & 0 & 0 & 0 \\
\rho_0        &-\alpha_0\bar\alpha_1&-\alpha_0\rho_1        & 0 & 0 & 0 & 0 & 0 \\
0 & \rho_1\bar\alpha_2 &-\alpha_1\bar\alpha_2 &\rho_2\bar\alpha_3& \rho_2\rho_3& 0 & 0 & 0 \\
0 & \rho_1\rho_2        &-\alpha_1\rho_2       &-\alpha_2\bar\alpha_3&-\alpha_2\rho_3 & 0 & 0 & 0 \\
0 & 0 & 0 & \rho_3\bar\alpha_4 &-\alpha_3\bar\alpha_4& \rho_4\bar\alpha_5 &\rho_4\rho_5 & 0 \\
0 & 0 & 0 &\rho_3\rho_4 & -\alpha_3\rho_4 &-\alpha_4\bar\alpha_5& -\alpha_4\rho_5 & 0 \\
0 & 0 & 0 & 0 & 0 & \rho_5\bar\alpha_6 &-\alpha_5\bar\alpha_6& \rho_6\bar\alpha_7 \\
0 & 0 & 0 & 0 & 0 & \rho_5\rho_6 & -\alpha_5\rho_6 &
-\alpha_6\bar\alpha_7
\end{bmatrix}
\end{equation*}
\caption{An $8\times 8$ CMV matrix in terms of the Verblunsky Coefficients}\label{F1}
\end{figure}

We will now give a brief overview of the contents of the paper.  We will structure this
presentation around the order of the sections.

In Section~2, we describe the origin of CMV matrices in the theory of orthogonal polynomials
on the unit circle.  Nothing said there is new; it is provided for completeness.

In Section~3, we show how a general unitary matrix can be reduced to a (direct sum) of
CMV matrices with a simple efficient algorithm.  Our model here is the famous Householder
implementation of the Lanczos reduction of a general symmetric matrix to tri-diagonal shape
by orthogonal conjugation \cite[\S 6.4]{Househ}.  This significantly reduces the storage
requirements during eigenvalue computations.

In Section~4, we describe a symplectic structure on the set of CMV matrices of fixed determinant.
Note that in the definition of CMV matrices, $\det(\Xi_k) = -1$ for $0\leq k \leq n-2$ and hence
$$
  \det\C=(-1)^{n-1}\bar\alpha_{n-1}.
$$
The specific symplectic structure is defined by algebraic means. This work is inspired by corresponding
work on Jacobi matrices stimulated by their appearance in the solution of the Toda lattice, \cite{Flaschka}.
Kostant, \cite{Kostant}, noticed that by altering the natural Lie bracket on real $n\times n$ matrices, the set of Jacobi
matrices with fixed trace becomes a co-adjoint orbit.  As a consequence, this manifold is a symplectic leaf
of the Lie--Poisson bracket.  Moreover, he uncovered an algebraic interpretation of the complete integrability
of the Toda lattice, which he extended to other finite dimensional Lie algebras.

As noted in Section~4, the natural algebraic setting for CMV matrices is a Lie group rather than
a Lie algebra.  We show that the set of CMV matrices with fixed determinant forms a symplectic leaf for
a certain Poisson structure on $GL(n;\Cmplx)$.  We will refer to this as the Gelfand--Dikij bracket,
though the name Sklyanin bracket would be equally appropriate.

In a concurrent paper, \cite{Li}, Luen-Chau Li has independently derived the main results of this
section.  In a sense, his approach is the reverse of ours: he studies the action of certain
dressing transformations, while we arrive at their existence only after studying the problem by
other means.

Section~4 also contains a few simple remarks regarding the Ablowitz-Ladik hierarchy. The
original (defocusing) Ablowitz Ladik equation \cite{AL1,AL2} is a space-discretization of the cubic nonlinear
Schr\"odinger equation:
\begin{equation}\label{ALeqn1}
-i \dot\beta_k = \rho_k^2 ( \beta_{k+1} + \beta_{k-1} ) - 2\beta_k,
\end{equation}
where $\beta_n$ is a sequence of numbers in the unit disk indexed over $\Ints$ and
$\rho_k=(1-|\beta_k|^2)^{1/2}$. If we change variables to $\alpha_k(t)=e^{2it}\beta_k(t)$, this system becomes
\begin{equation}\label{ALeqn}
-i \dot\alpha_k = \rho_k^2 ( \alpha_{k+1} + \alpha_{k-1} ),
\end{equation}
which is a little simpler.  If we then choose $\alpha_{-1}$ and $\alpha_{n-1}$ to lie
on the unit circle, then they do not move and we obtain a finite system of ODEs for
$\alpha_k$, $0\leq k \leq n-2$, which is the specific case we treat.

These equations form a completely integrable Hamiltonian system; indeed they were introduced
with this very property in mind.  By considering all possible functions of the commuting
Hamiltonians, one is immediately lead to a hierarchy of equations containing \eqref{ALeqn}
as a special case.  Similar considerations lead from the original Toda equation to the full
Toda hierarchy.

As presented in Section~4, the matters discussed under the rubric `Ablowitz--Ladik hierarchy'
are not obviously linked to the equations just discussed.  We make this connection in Section~5.

As a Hamiltonian system, the Ablowitz--Ladik equation comes with a symplectic form.  In Section~5,
we prove that this agrees with the bracket we introduced from algebraic considerations.  Viewed from
the opposite perspective, the goal of this section is to write the Gelfand--Dikij bracket tensorially
using the Verblunsky coefficients as a system of coordinates.

In Sections~\ref{sESM} and~\ref{sAS}, we study further properties of the Ablowitz--Ladik system
following the work of Moser, \cite{Moser}, on the Toda lattice.  Specifically, we study how
the spectral measure naturally associated to a CMV matrix evolves under the Hamiltonians
of the Ablowitz--Ladik hierarchy.  This information is then used to derive long-time asymtotics
and determine the scattering map.

The last topic we treat is the construction of action/angle coordinates for the finite Ablowitz--Ladik
system.  This is done in Section~8.

\medskip
\noindent\textit{Acknowledgements:} We are delighted to dedicate this paper to Percy Deift
on the occasion of his sixtieth birthday.  Our understanding of the matters discussed here
owes much to
his lectures, \cite{Deift}.  We are also grateful to him for encouragement
along the way.

\section{Orthogonal Polynomials.}\label{s2}

As CMV matrices arose in the study of orthogonal polynomials, it is
natural that we begin there. We will first describe the relation of
orthogonal polynomials to Jacobi matrices and then explain the
connection to CMV matrices.

Given a probability measure $d\nu$ supported on a finite subset of
$\Reals$, say of cardinality $n$, we can apply the Gram--Schmidt
procedure to $\{1,x,x^2,\ldots,x^{n-1}\}$ and so obtain an
orthonormal basis for $L^2(d\nu)$ consisting of
polynomials, $\{p_j(x):j=0,\ldots,{n-1}\}$, with positive leading
coefficient. In this basis, the linear transformation $f(x)\mapsto
xf(x)$ is represented by a Jacobi matrix.  An equivalent statement
is that the orthonormal polynomials obey a three-term recurrence:
$$
xp_{j}(x) = a_{j} p_{j+1}(x) + b_{j}p_j(x) + a_{j-1} p_{j-1}(x)
$$
where $a_{-1}=0$ and $p_{n}\equiv 0$.  A third equivalent statement
is the following:  $\lambda$ is an eigenvalue of $J$ if and only if
$\lambda\in\supp(d\nu)$; moreover, the corresponding eigenvector is
$[p_{0}(\lambda),p_{1}(\lambda),\ldots,p_{n-1}(\lambda)]^{T}$.

We have just shown how measures on $\Reals$ lead to Jacobi matrices;
in fact, there is a one-to-one correspondence between them.  Given a
Jacobi matrix, $J$, let $d\nu$ be the spectral measure associated to
$J$ and the vector $e_1=[1,0,\ldots,0]^{T}$.  Then $J$ represents
$x\mapsto xf(x)$ in the basis of orthonormal polynomials associated
to $d\nu$.

Before explaining the origin of CMV matrices, it is necessary to
delve a little into the theory of orthogonal polynomials on the unit
circle.  For a more complete description of what follows, the reader
should turn to \cite{Simon1}.

Given a finitely-supported probability measure $d\mu$ on $S^1$, the
unit circle in $\Cmplx$, we can construct an orthonormal system of
polynomials, $\phi_k$, by applying the Gram--Schmidt procedure to
$\{1,z,\ldots\}$.  These obey a recurrence relation; however, to
simplify the formulae, we will present the relation for the monic
orthogonal polynomials $\Phi_k(z)$:
\begin{align}
\Phi_{k+1}(z)   &= z\Phi_k(z)   - \bar\alpha_k \Phi_k^*(z).
\label{PhiRec}
\end{align}
Here $\alpha_k$ are recurrence coefficients, which are called
Verblunsky coefficients, and $\Phi_k^*$ denotes the reversed
polynomial:
\begin{equation}\label{rev}
\Phi_k(z) = \sum_{l=0}^k c_l z^l \quad \Rightarrow \quad \Phi_k^*(z)
= \sum_{l=0}^k \bar{c}_{k-l} z^l.
\end{equation}
When $d\mu$ is supported at exactly $n$ points, $\alpha_k\in\Disk$
for $0\leq k\leq {n-2}$ while $\alpha_{n-1}$ is a unimodular complex
number.  (Incidentally, if $d\mu$ has infinite support, then there are
infinitely many Verblunsky coefficients and all lie inside the unit disk.)

To recover the relation between the orthonormal polynomials, one
need only apply the following relation, which can be deduced from
\eqref{PhiRec}:
\begin{equation}\label{E:PhiNorm}
\bigl\| \Phi_k \bigr\|_{L^2(d\mu)} = \prod_{l=0}^{k-1} \rho_l,
    \qquad\text{where $\rho_l=\sqrt{1 - |\alpha_l|^2}$.}
\end{equation}

The Verblunsky coefficients completely describe the measure $d\mu$:

\begin{theorem}[Verblunsky]\label{T:Verbl}
There is a 1-to-1 correspondence between probability measures on the unit circle
supported at $n$ points and Verblunsky coefficients $(\alpha_0,\ldots,\alpha_{n-1})$
with $\alpha_{k}\in\ID$ for $0\leq k\leq n-2$ and $\alpha_{n-1}\in S^1$.
\end{theorem}

From the discussion of Jacobi matrices, it would be natural to consider the matrix representation
of $f(z)\mapsto zf(z)$ in $L^2(d\mu)$ with respect to the basis of orthonormal polynomials. This is
\textit{not} a CMV matrix; rather it is what Simon, \cite{Simon1}, has dubbed a GGT matrix, from the
initials of Geronimus, Gragg, and Teplyaev.  Perhaps the most striking difference from a CMV (or
Jacobi) matrix is that a GGT matrix is very far from sparse---generically, all entries above and
including the sub-diagonal are non-zero. (In a sense, CMV matrices are optimally sparse; see
Theorem~\ref{Tsparse}.)

Cantero, Moral, and Velazquez had the simple and ingenious idea of applying the Gram--Schmidt
procedure to $\{1,z,z^{-1},z^2,z^{-2},\ldots\}$ rather than $\{1,z,\ldots\}$. The resulting
functions, $\chi_k(z)$ ($0\leq k\leq n-1$), are easily expressed in terms of the orthonormal polynomials:
\begin{equation}
\chi_k(z) = \begin{cases} z^{-k/2} \phi_k^*(z) & \text{: $k$ even} \\
        z^{-(k-1)/2} \phi_k(z) & \text{: $k$ odd.} \end{cases}
\end{equation}
In this basis, the map $f(z)\mapsto zf(z)$ is represented in an especially simple form:

\begin{theorem}\label{T:2a}
In the orthonormal basis $\{\chi_k(z)\}$ of $L^2(d\mu)$, the operator $f(z)\mapsto zf(z)$ is
represented by the CMV matrix associated to the Verblunsky coefficients of the measure $d\mu$.
\end{theorem}

The $\C=\mathcal{LM}$ factorization presented in the introduction originates as follows:
Let us write $x_k$, $0\leq k\leq n-1$, for the orthonormal basis constructed by applying the Gram--Schmidt
procedure to $\{1,z^{-1},z,z^{-2},z^2,\ldots\}$.  Then the matrix elements of $\mathcal{L}$ and
$\mathcal{M}$ are given by
$$
\mathcal{L}_{j+1,k+1}=\langle \chi_j(z) |z  x_k(z) \rangle, \qquad
\mathcal{M}_{j+1,k+1}=\langle x_j(z) | \chi_k(z) \rangle.
$$
See \cite{Simon1} for further discussion.

The measure $d\mu$ can be reconstructed from $\C$ in a manner analogous to the Jacobi case:

\begin{theorem}\label{T:2b}
Let $d\mu$ be the spectral measure associated to a CMV matrix, $\C$,
and the vector $e_1$.  Then $\C$ is the CMV matrix associated to the
measure $d\mu$.
\end{theorem}

Proofs of these two Theorems can be found in
\cite{CMV} or \cite{Simon1}. The second result also
follows from Corollary~\ref{C:DUCK} below.

\section{Reduction to CMV shape.}\label{s3}

The purpose of this section is to explain a simple analytic (and
numeric) algorithm for the reduction of a unitary matrix to CMV
form. At the core of the algorithm is the Householder method for
reducing the number of nonzero entries of a matrix by conjugating it
with reflections. These are chosen in such a way
as to allow for successive applications of the reduction.

\begin{figure}[h]
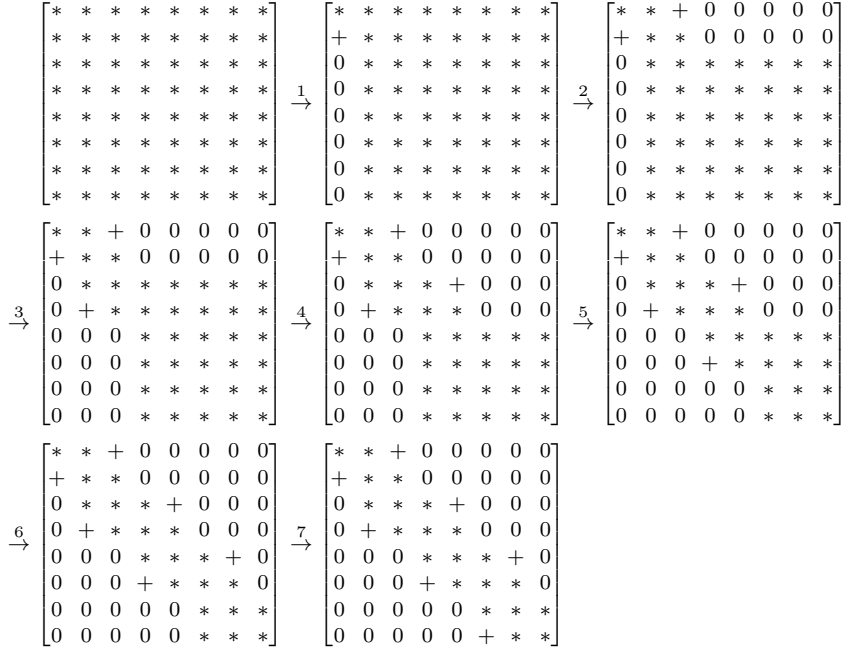

\def\p{\hbox to 0mm{\hss+\hss}}
{ \footnotesize \setlength{\arraycolsep}{1ex}
\begin{alignat*}{1}
\begin{bmatrix}
* & * & * & * & * & * & * & * \\ * & * & * & * & * & * & * & * \\
* & * & * & * & * & * & * & * \\ * & * & * & * & * & * & * & * \\
* & * & * & * & * & * & * & * \\ * & * & * & * & * & * & * & * \\
* & * & * & * & * & * & * & * \\ * & * & * & * & * & * & * & *
\end{bmatrix} &\overset{1}\to \begin{bmatrix}
* & * & * & * & * & * & * & * \\\p & * & * & * & * & * & * & * \\
0 & * & * & * & * & * & * & * \\ 0 & * & * & * & * & * & * & * \\
0 & * & * & * & * & * & * & * \\ 0 & * & * & * & * & * & * & * \\
0 & * & * & * & * & * & * & * \\ 0 & * & * & * & * & * & * & *
\end{bmatrix} \overset{2}\to \begin{bmatrix}
* & * & \p& 0 & 0 & 0 & 0 & 0 \\ \p& * & * & 0 & 0 & 0 & 0 & 0 \\
0 & * & * & * & * & * & * & * \\ 0 & * & * & * & * & * & * & * \\
0 & * & * & * & * & * & * & * \\ 0 & * & * & * & * & * & * & * \\
0 & * & * & * & * & * & * & * \\ 0 & * & * & * & * & * & * & *
\end{bmatrix} \\
\overset{3}\to \begin{bmatrix}
* & * & \p& 0 & 0 & 0 & 0 & 0 \\ \p& * & * & 0 & 0 & 0 & 0 & 0 \\
0 & * & * & * & * & * & * & * \\ 0 & \p& * & * & * & * & * & * \\
0 & 0 & 0 & * & * & * & * & * \\ 0 & 0 & 0 & * & * & * & * & * \\
0 & 0 & 0 & * & * & * & * & * \\ 0 & 0 & 0 & * & * & * & * & *
\end{bmatrix} &\overset{4}\to \begin{bmatrix}
* & * & \p& 0 & 0 & 0 & 0 & 0 \\ \p& * & * & 0 & 0 & 0 & 0 & 0 \\
0 & * & * & * & \p& 0 & 0 & 0 \\ 0 & \p& * & * & * & 0 & 0 & 0 \\
0 & 0 & 0 & * & * & * & * & * \\ 0 & 0 & 0 & * & * & * & * & * \\
0 & 0 & 0 & * & * & * & * & * \\ 0 & 0 & 0 & * & * & * & * & *
\end{bmatrix} \overset{5}\to \begin{bmatrix}
* & * & \p& 0 & 0 & 0 & 0 & 0 \\ \p& * & * & 0 & 0 & 0 & 0 & 0 \\
0 & * & * & * & \p& 0 & 0 & 0 \\ 0 & \p& * & * & * & 0 & 0 & 0 \\
0 & 0 & 0 & * & * & * & * & * \\ 0 & 0 & 0 & \p& * & * & * & * \\
0 & 0 & 0 & 0 & 0 & * & * & * \\ 0 & 0 & 0 & 0 & 0 & * & * & *
\end{bmatrix} \\
\overset{6}\to \begin{bmatrix}
* & * & \p& 0 & 0 & 0 & 0 & 0 \\ \p& * & * & 0 & 0 & 0 & 0 & 0 \\
0 & * & * & * & \p& 0 & 0 & 0 \\ 0 & \p& * & * & * & 0 & 0 & 0 \\
0 & 0 & 0 & * & * & * & \p& 0 \\ 0 & 0 & 0 & \p& * & * & * & 0 \\
0 & 0 & 0 & 0 & 0 & * & * & * \\ 0 & 0 & 0 & 0 & 0 & * & * & *
\end{bmatrix} &\overset{7}\to \begin{bmatrix}
* & * & \p& 0 & 0 & 0 & 0 & 0 \\ \p& * & * & 0 & 0 & 0 & 0 & 0 \\
0 & * & * & * & \p& 0 & 0 & 0 \\ 0 & \p& * & * & * & 0 & 0 & 0 \\
0 & 0 & 0 & * & * & * & \p& 0 \\ 0 & 0 & 0 & \p& * & * & * & 0 \\
0 & 0 & 0 & 0 & 0 & * & * & * \\ 0 & 0 & 0 & 0 & 0 & \p& * & *
\end{bmatrix}
\end{alignat*}
\vskip -1ex
} \caption{The CMV-ification algorithm in action.}\label{F2}
\end{figure}

The Householder algorithm is the method of choice to reduce any real symmetric
matrix to tridiagonal form or any non-symmetric matrix to Hessenberg shape.
Here we implement the successive reductions differently, by
alternating columns and rows (see Figure~\ref{F2}) to reduce any unitary
matrix for which $e_1=[1,0,\ldots,0]^T$ is cyclic to CMV shape.  With the obvious modification,
it reduces other unitary matrices to a direct sum of CMV matrices.

Let us briefly describe the construction of Householder reflections; what we present is
a very slight modification of what one will find in most textbooks.
Given $u\in\Cmplx^n$ and $0\leq m < n$, let
$$
v=\big[0,\ldots,0,\alpha,u_{m+2},\ldots,u_n\big]^T \ \ \text{with}\ \
\alpha=u_{m+1}-\frac{u_{m+1}}{|u_{m+1}|} \biggl(\sum_{j=m+1}^n
|u_j|^2\biggr)^{1/2}.
$$
The reflection through the plane perpendicular to $v$ is given by
$$
R=I-2\frac{vv^{\dagger}}{\|v\|^2}.
$$
Naturally, $R^{\dagger}=R^{-1}=R$; moreover,
$$
Ru=u-v=\big[u_1,\ldots,u_{m},u_{m+1}-\alpha,0\ldots,0\big]^T
$$
and since the first $m$ entries of $v$ are identically zero, $R$ has the block structure
\begin{equation}\label{Rstruct}
R=
\begin{bmatrix}
I_{m} & 0\\
0 & * \\
\end{bmatrix},
\end{equation}
where $I_{m}$ is the $m\times m$ identity matrix, and $*$ denotes an
unspecified $(n-m)\times(n-m)$ matrix.

When working over $\Reals$, one typically chooses $\alpha$ slightly differently to guarantee that
$u_{m+1}-\alpha$ is positive (or zero).  This is not possible in our case.  To remedy the situation,
we multiply by a suitably chosen rotation $D=\diag(1,\ldots,1,e^{i\theta},1,\ldots,1)$ where $e^{i\theta}$
sits in the $(m+1)$th position.

Replacing $R$ by $DR$ gives us the desired modification of the standard Householder
reflection.  We summarize its properties as follows:

\begin{lemma}
Given $u\in\Cmplx^n$ and $0\leq m<n$ there is a Householder reflection $R$ that has the
block structure \eqref{Rstruct} and
$$
R: u \mapsto \Bigl[u_1,\ldots,u_m, \sqrt{\textstyle\sum_{l=m+1}^n |u_l|^2 },0,\ldots,0\Bigr]^{T}.
$$
We will refer to $R$ as the reflection at level $m$ for $u$.
\end{lemma}

\begin{theorem}[A Householder algorithm for unitary matrices]\label{ReductionCMV}
Any $n\times n$ unitary matrix, $U$, for which $e_1$ is cyclic can be reduced to
CMV shape by successive conjugations with Householder reflections.  Moreover, this
can be done without altering the spectral measure associated to $e_1$.
\end{theorem}

\begin{proof}
During the proof, we will refer the reader to Figure~\ref{F2}.  The small numbers above each arrow
enumerate the steps in the algorithm.

As we will see, at every step $U$ is conjugated by a matrix that leaves $e_1$ invariant. Thus
the spectral measure associated to $e_1$ will not change.  In particular, this vector
will remain cyclic throughout.

For the first step, we choose $R$ to be the reflector at level two for the first column of $U$
and then form $RU\!R^\dagger$.  Left multiplication by $R$ creates zeros in the places indicated
and these are not destroyed by right multiplication by $R^\dagger$ because of its block structure.
One detail remains: why the entry marked $+$ is not zero.  The answer is simple.  This entry would
only vanish if the first column of the original matrix where a multiple of $e_1$.  This
would make $e_1$ an eigenvector for $U$ and so not cyclic.

Next we describe Step~3 as an archetype for all subsequent odd-numbered steps.  Even-numbered steps
differ from odd-numbered steps only by interchanging the roles of rows and columns.  We will write
$U$ for the matrix output by step two (not the original one).

Let $R$ be the reflector at level 3 associated to the second column of $U$, then the effect of
Step~3 is $U\mapsto RU \!R^\dagger$. For the reasons explained in Step~1, this will produce zeros in the
second column as shown in Figures~\ref{F2} and~\ref{F3}. This state of knowledge is represented by the
centre matrix in Figure~\ref{F3}.  Two questions remain however:
why is the entry newly marked $+$ not zero, and why do additional zeros appear in the fourth column.
We answer them in reverse order.

\begin{figure}[h]
\def\p{\hbox to 0mm{\hss+\hss}}
{ \footnotesize \setlength{\arraycolsep}{1ex}
\begin{alignat*}{1}
\begin{bmatrix}
* & * & \p& 0 & 0 & 0 & 0 & 0 \\ \p& * & * & 0 & 0 & 0 & 0 & 0 \\
0 & * & * & * & * & * & * & * \\ 0 & * & * & * & * & * & * & * \\
0 & * & * & * & * & * & * & * \\ 0 & * & * & * & * & * & * & * \\
0 & * & * & * & * & * & * & * \\ 0 & * & * & * & * & * & * & *
\end{bmatrix} &\overset{3}\to \begin{bmatrix}
* & * & \p& 0 & 0 & 0 & 0 & 0 \\ \p& * & * & 0 & 0 & 0 & 0 & 0 \\
0 & * & * & * & * & * & * & * \\ 0 & + & * & * & * & * & * & * \\
0 & 0 & * & * & * & * & * & * \\ 0 & 0 & * & * & * & * & * & * \\
0 & 0 & * & * & * & * & * & * \\ 0 & 0 & * & * & * & * & * & *
\end{bmatrix} \overset{?}= \begin{bmatrix}
* & * & \p& 0 & 0 & 0 & 0 & 0 \\ \p& * & * & 0 & 0 & 0 & 0 & 0 \\
0 & * & * & * & * & * & * & * \\ 0 & \p& * & * & * & * & * & * \\
0 & 0 & 0 & * & * & * & * & * \\ 0 & 0 & 0 & * & * & * & * & * \\
0 & 0 & 0 & * & * & * & * & * \\ 0 & 0 & 0 & * & * & * & * & *
\end{bmatrix}
\end{alignat*}
\vskip -1ex } \caption{Matrices in the discussion of Step~3.}\label{F3}
\end{figure}

The rows of a unitary matrix must be orthogonal; in particular, all rows must be perpendicular
to the first.  As the third entry in the top row is non-zero (indeed positive), it follows that
the third column must have zeros in the places indicated (cf. the third matrix in Figure~\ref{F3}).

If the entry marked $+$ in the second column were actually zero, then by the reasoning of the
previous paragraph, its right-hand neighbour would also vanish.  This contradicts the cyclicity
of $e_1$; specifically, the linear span of $\{e_1,e_2,e_3\}$ would be an invariant subspace.

With obvious modifications, these arguments apply to any odd step ($\geq 3$).
As noted even steps are analogous but with rows and columns reversed; indeed
one can just apply the discussion above to $U^\dagger$.  (As $U$
is unitary, $e_1$ is also cyclic for $U^\dagger$.)

Note that in the last step, no additional zeros are produced.  This does not result
in any modifications to the argument, we merely wish to point out that in this case, $R$ is a diagonal
matrix.
\end{proof}

Theorem~\ref{ReductionCMV} tells us how to reduce any
unitary matrix having $e_1$ as a cyclic vector to CMV shape. Now we wish to show
that any unitary matrix $U$ in CMV shape is actually a CMV matrix,
i.e., there exists a set of coefficients
$\alpha_0,\ldots,\alpha_{n-2}\in\ID$ and $\alpha_{n-1}\in S^1$ so
that $U$ is the matrix given by \eqref{F:CMV}. This result is not new,
and can be found in \cite{CMV2}. We prove this result differently, and
obtain it as a corollary of the following

\begin{prop}\label{UniqCMV}
Let $B$ and $C$ be two unitary matrices in CMV shape that have the
same spectral measure with respect to the vector $e_1$. Then $B=C$.
\end{prop}

\begin{proof}
Since $B$ and $C$ have the same spectral measure with respect to
$e_1$, there exists a unitary matrix $V$ so that $Ve_1=e_1$ and
\begin{equation}\label{conj}
C=V^\dagger B V.
\end{equation}
We inductively prove that $Ve_k=e_k$ for all $1\leq k\leq n$, which
shows that $V$ is the identity and so $B=C$.

Assume that $Ve_k=e_k$ holds for all $k\leq m$.  (From the
definition of $V$, it is true when $m=1$.) This implies that the top
left $m\times m$ blocks of $B$ and $C$ coincide:
$$
\langle e_j\vert C e_k\rangle =\langle Ve_j\vert B (Ve_k)\rangle
    = \langle e_j\vert B e_k\rangle
$$
for all $1\leq j,k \leq m$.  We now proceed differently depending on
$m$.

If $m\geq 3$ is odd, consider the $(m-1)^{\text{th}}$ columns of $C$ and $B$.
As both matrices are of CMV shape and agree on the principle
$m\times m$ minor, the only place that these columns may differ is
in the $(m+1)^\text{th}$ position.  In fact, they must agree here too
because CMV shape requires that this entry be positive and then it
is uniquely determined by the fact that the columns of a unitary
matrix have norm one.

From $VC=BV$, $Ve_{m-1}=e_{m-1}$, and $C e_{m-1}=Be_{m-1}$, which we
just proved, we have
\begin{equation}\label{E:col}
(V-1)C e_{m-1} = (BV-C)e_{m-1} = 0.
\end{equation}
But now by the shape of the $(m-1)^\text{th}$ column of $C$ described above,
we know that
$$
C e_{m-1} = C_{m+1,m-1} e_{m+1} + \sum_{k=m-2}^{m} C_{k,m-1} e_{k}
$$
with $C_{m+1,m-1}>0$.  Substituting this into \eqref{E:col} gives
$(1-V)e_{m+1}=0$, which concludes the inductive step in this case.

If $m=1$, one should repeat the above analysis using the first
column.

If $m$ is even, then one must repeat the above argument using rows
instead of columns. Alternatively, one may proceed as above after
first taking the (conjugate) transpose in \eqref{conj}.
\end{proof}

\begin{coro}\label{C:DUCK}
Let $U$ be an $n\times n$ unitary matrix in CMV shape.
Then there is a unique $(n-1)$-tuple of Verblunsky coefficients
$\alpha_0,\ldots,\alpha_{n-2}\in\ID$ and $\alpha_{n-1}\in S^1$ so
that $U$ is the CMV matrix associated to these coefficients.
\end{coro}

\begin{proof}
For a unitary matrix in CMV shape, let $\mu$ be the spectral measure
associated to $e_1$, and $\alpha_0,\ldots,\alpha_{n-1}$ the
corresponding Verblunsky coefficients. The CMV matrix $\C$ given by
Definition~\ref{D:CMV} has $\mu$ as its spectral measure. So, by
Proposition~\ref{UniqCMV}, the original matrix must equal $\C$,
which is indeed a CMV matrix. Uniqueness of the coefficients
is a direct consequence of Proposition~\ref{UniqCMV} and
Verblunsky's Theorem~\ref{T:Verbl}.
\end{proof}

\begin{remark}
There is another way of proving the corollary, which involves
recursively identifying the $\alpha$'s in a unitary matrix of CMV
shape. This is essentially the proof of Theorem~3.8 from
\cite{CMV2}.

Further note that one can read the Verblunsky coefficients from a
CMV matrix. In other words, for any unitary matrix $C$ in CMV shape
there exists a \emph{unique} set of Verblunsky coefficients so
that $C$ is the CMV matrix with these coefficients.
\end{remark}

We close this section by presenting a result of Cantero, Moral, and Velazquez
(see Theorem~3.9 of \cite{CMV2}) which implies that, in terms of sparseness, CMV matrices are minimal among
unitary matrices having cyclic vectors.

We call a unitary matrix, $U$, $(p,q)$-diagonal if it has only $p$ nonzero subdiagonals and $q$ nonzero
superdiagonals. In other words, we require that $U_{j,k}=0$ for $j-k\geq p+1$ and $k-j\geq q+1$. Note that
CMV matrices are $(2,2)$-diagonal. With this definition, the statement is the following:

\begin{theorem}\label{Tsparse}
A unitary $(p,1)$-diagonal matrix is a sum of diagonal blocks of
order no greater than $p+1$.  The same is true of any unitary $(1,p)$-diagonal matrix.
\end{theorem}

In particular, note that any 4-diagonal unitary matrix is a direct sum of $3\times3$ blocks.  As breaking
a matrix into a direct sum is tantamount to factoring the characteristic polynomial, no exact algorithm
can perform this task for generic input.  Similarly, it is impossible to reduce a unitary
matrix to four diagonals whilst preserving the spectral measure associated to $e_1$.

\section{Lie structure.}\label{s4}

The manifold of Jacobi matrices with fixed trace forms a co-adjoint
orbit associated to a particular Lie algebra structure on the
$n\times n$ matrices (although not the one induced by
$GL(n;\Reals)$).  This gives the manifold a symplectic structure.
These matters are described in detail in
\cite{Deift,IntSysII,Perel}, for example.

In this section, we will show that the manifold of CMV matrices with
fixed determinant is a symplectic leaf for a natural Poisson
structure on $GL(n,\Cmplx)$ and (what is basically equivalent) it is
exactly the orbit of a natural group action on this space. and so
also have a natural symplectic structure.  For pedagogical reasons,
we will present the symplectic structure first and then describe how
it relates to a group action.  (As described in the Introduction,
\cite{Li}, takes the opposite approach to arrive at the same conclusions.)

In deference to those new to the subject, we will endeavour to use
the notation and terminology of \cite{IntSysII}.  We will also
mostly refer to this book for proofs that we omit; we hope that the
authors of the original articles will forgive this indirect
reference to their work.

Self-adjoint matrices form a vector space, hence it is natural that the
orbit and symplectic structures on Jacobi matrices arises from the study
of a Lie algebra (actually, its vector-space dual). Unitary matrices
form a group and hence the natural backdrop for CMV is a Lie group.
Specifically, the construction we give below is that of a Sklyanin
bracket on $Gl(n,\Cmplx)$; see \cite[\S 2.12]{IntSysII}. However, we
choose to give a presentation in which the algebra of matrices takes
centre stage; an analogous construction for KdV using the algebra of
pseudo-differential operators was given by Gelfand and Dickij
\cite{GD}. This approach is described in Section~2.12.6 of
\cite{IntSysII}.  (Note that here we are referring to the second
symplectic structure associated with KdV, which was originally
proposed by Adler \cite[\S 4]{Adler}.)

Let $\g$ denote the (associative) algebra of $n\times n$ complex
matrices.  The algebra structure gives rise to a natural Lie algebra
structure:
$$
  [B,C] = BC-CB.
$$
Of course, this also results from viewing $\g$ as the Lie algebra of
$GL(n;\Cmplx)$.

As a vector space, $\g=\LL\oplus\A$, where
$$
\A = \{A : A = - A^\dagger \},
$$
which is the Lie algebra of the group $\U(n)$ of $n\times n$ unitary
matrices, and
$$
\LL = \{A\in\g : L_{i,j} = 0 \text{ for $i>j$ and } L_{i,i}\in\Reals
\}
$$
which is the Lie algebra of the group $\LLL(n)$ of $n\times n$ lower
triangular matrices with positive diagonal entries. We will write
$\pi_{\A}$ and $\pi_\LL$ for the projections into these summands.

This vector-space splitting of $\g$ permits us to give it a second
Lie algebra structure. First we define $R:\g\to\g$ by either
\begin{equation}\label{E:RDef}
\begin{aligned}
R(X) &= \pi_\LL(X) - \pi_\A(X),  &  \qquad&\text{for all $X\in\g$, or} \\
R(L+A)&= L-A    &  \qquad&\text{for all $A\in\A$ and $L\in\LL$.}
\end{aligned}
\end{equation}
The second Lie bracket can then be written as either
\begin{equation}\label{E:DefnRBracket}
\begin{aligned}{}
[X,Y]_{\RR} = \tfrac12 [R(X),Y] + \tfrac12 [X,R(Y)]     \quad&\text{$\forall$ $X,Y\in\g$, or} \\
[L+A,L'+A']_{\RR} = [L,L'] - [A,A']    \quad&\text{$\forall$
$L,L'\in\LL$, and $A,A'\in\A$.}
\end{aligned}
\end{equation}
The second definition also makes it transparent that the $R$-bracket
obeys the Jacobi identity.

Whether one treats Jacobi or CMV matrices, the (vector space) dual
$\g^*$ of $\g$ plays an important role. We can identify it with $\g$
using the pairing
\begin{equation}\label{pairing}
\langle X,Y\rangle = \Im \tr (XY).
\end{equation}
This is not an inner product; however, it is non-degenerate---just
choose $Y=iX^\dagger$.  It is also symmetric and ad-invariant:
\begin{equation}
\langle X, [Z,Y] \rangle = \Im \tr (XZY-XYZ) = \langle
[X,Z],Y\rangle.
\end{equation}
This is equivalent to Ad-invariance:
\begin{equation}
\langle BXB^{-1}, BYB^{-1} \rangle = \langle X, Y
\rangle,\qquad\text{for any $B\in GL(n,\Cmplx)$.}
\end{equation}

Concomitant with the pairing is an identification of 1-forms and
vector fields on~$\g$:

\begin{definition}
Given $\phi:\g\to\Reals$ and $B\in\g$, define $\nabla\phi:\g\to\g$
by
\begin{equation}
\frac{d}{dt}\bigg|_{t=0} \phi(B+tC) =
\langle\nabla\phi\big|_B,C\rangle.
\end{equation}
Equivalently, if we write $b_{k,l}=u_{k,l} + i v_{k,l}$ for the
matrix entries of $B$, then
\begin{equation}
 [\nabla\phi]_{k,l} = \frac{\partial\phi}{\partial v_{l,k}} + i \frac{\partial\phi}{\partial u_{l,k} }
\end{equation}
(notice the reversal of the order of the indices).
\end{definition}

We can now define the desired Poisson bracket on $\g$. In fact, by
first checking some elementary properties of the objects set out
above, one can show that it obeys the Jacobi identity.

\begin{prop}\label{P:PS}
Let $R:\g\to\g$ and $\langle\cdot\,,\cdot\rangle$ be as above.
\begin{SL}
\item $\A^\perp=\A$ and $\LL^\perp=\LL$.
\item $R$ is antisymmetric: $\langle X,R(Y) \rangle = - \langle R(X),Y \rangle$.
\item $R$ obeys the modified classical Yang-Baxter equation:
\begin{equation*}
[R(X),R(Y)] - R\bigl([R(X),Y] + [X,R(Y)]\bigr) = -[X,Y]
\end{equation*}
\item Given $\phi,\psi:\g\to\Reals$, let $X=\nabla\phi$ and $Y=\nabla\psi$.  Then
\begin{align}\label{E:GDD}
\{\phi,\psi\}\bigr|_B &= \tfrac{1}{2} \langle R(XB), YB \rangle - \tfrac{1}{2}\langle R(BX), BY \rangle \\
&=\tfrac{1}{2} \langle R(XB), YB \rangle + \tfrac{1}{2}\langle BX,
R(BY) \rangle.
\end{align}
defines a Poisson structure on $\g$.
\end{SL}
\end{prop}

\begin{proof} (i) For $A,B\in \A$,
\begin{align*}
\Im\tr(AB)=\tfrac1{2i}\tr(AB-B^\dagger A^\dagger) =
\tfrac1{2i}\tr(AB-BA) = 0.
\end{align*}
For $L,L'\in\LL$, the diagonal entries of $LL'$ are products of the
corresponding entries in $L$ and $L'$. In particular, they are real;
therefore, $\Im \tr(LL')=0$.

We have just shown that $\A^\perp\subseteq\A$ and
$\LL^\perp\subseteq\LL$.  As the pairing is non-degenerate, equality
follows by dimension counting:
$\dim(\A)=\dim(\LL)=n^2=\tfrac12\dim(\g)$.

(ii) Antisymmetry follows from part (i) by simple computation:
\begin{align*}
\langle R(L+A),L'+A'\rangle = \langle L,L'\rangle - \langle
A,A'\rangle = \langle L+A, R(L'+A')\rangle.
\end{align*}

(iii)  This can be checked by direct computation.  However, this
belies the important role of the mCYB equation in the theory of
double Lie algebras and thence in the theory of integrable systems.
For more information, see \cite[\S 2.2]{IntSysII} or
\cite{Deift,Perel}.

(iv) This is an example of a much more general statement, namely,
these formulae define a Poisson bracket on $\g$ whenever
$\langle\cdot\,,\cdot\rangle$ is symmetric and ad-invariant and $R$
has properties (ii) and (iii). The only axiom of a Poisson bracket
that is not immediate from the definition is the Jacobi identity.
Unfortunately, checking this is consumes more space than we can
justify here.  Details can be found in \cite[\S 4]{GD}, but readers
should not be discouraged from simply doing it themselves. We hope
that the following two hints make this more palatable:
\begin{gather}
\nabla\{\phi,\psi\}\bigl|_B = \tfrac12 R(XB)Y - \tfrac12 R(YB)X +
\tfrac12 XR(BY) - \tfrac12 YR(BX)
\end{gather}
and given any trio of matrices $A$, $B$, and $C$,
\begin{align*}
&\hphantom{{}={}} \sum \langle R(B)A-R(A)B,\ R(C) \rangle \\
&= \sum \tfrac23 \langle R\bigl(R(B)A-R(A)B\bigr) ,\  C \rangle
    +\tfrac13 \langle R(B)R(A) -  R(A)R(B),\ C \rangle \\
&=  \tfrac13 \sum  \langle [A,B],\ C \rangle
\end{align*}
where all sums are over cyclic permutations of $(A,B,C)$.
\end{proof}

The Poisson bracket given in the proposition is often referred to as
the `quadratic bracket' because the point $B\in\g$ where it is
evaluated appears quadratically. The Lie-Poisson (or Kirillov)
bracket associated to the Lie algebra $(\g,[\,,]_\RR)$ is linear
in~$B$:
$$
\{\phi,\psi\}_{{}_{LP}}\bigr|_B = \tfrac{1}{2} \langle [X,Y]_{\RR} ,
B \rangle.
$$
This bracket is known to be compatible (in the sense of
Magri-Lenard) with the bracket above; however, it is not pertinent
to the study of CMV matrices.  It is relevant to Jacobi matrices:
under the embedding $J\mapsto iJ$, the manifolds of Jacobi matrices
with fixed trace are symplectic leaves. This is just the usual
construction in $Gl(n,\Reals)$ in disguise.

\begin{lemma}\label{L:HV}
The Hamiltonian vector field on $\g$ associated to
$\phi:\g\to\Reals$ is
\begin{equation}\label{E:HV}
\dot B = \tfrac12 \bigl[ BR(XB) - R(BX)B \bigr]
\end{equation}
where $X=\nabla \phi$.  Equivalently,
\begin{equation}\label{E:HV2a}
\dot B = B\pi_\LL(XB) - \pi_\LL(BX)B = \pi_\A(BX)B - B\pi_\A(XB) .
\end{equation}
In particular, if $B$ is unitary, then
\begin{equation}\label{E:HV-U}
\dot B = - B \,\pi_\A(B^{-1}LB)
\end{equation}
where $L=\pi_\LL(BX)$.
\end{lemma}

\begin{proof}
In terms of the Poisson bracket, the defining property of the
Hamiltonian vector field is $\{\phi,\psi\} = \dot \psi = \langle
\dot B, \nabla\psi\rangle$.  Equation \eqref{E:HV} now follows by
cycling the trace: if $Y=\nabla\psi$, then
\begin{align}
\{\phi,\psi\} &= \tfrac{1}{2} \langle R(XB), YB \rangle - \tfrac{1}{2}\langle R(BX), BY \rangle  \\
&= \tfrac{1}{2}\langle BR(XB) - R(BX)B , Y \rangle.
\end{align}

From the definition, $R=\Id-2\pi_\A=2\pi_\LL-\Id$.  Equation
\eqref{E:HV2a} follows by substituting these relations into
\eqref{E:HV}.

To obtain the special case, let us write $BX=L+A$. As $A$ and
$B^{-1}AB$ are anti-Hermitian, $\pi_\A(A)B = B\pi_\A(B^{-1}AB)$.
This implies
\begin{equation}
\dot B = \pi_\A\bigl(L+A\bigr)B - B\pi_\A\bigl(B^{-1}(L+A)B\bigr) =
- B\pi_\A\bigl(B^{-1}LB\bigr),
\end{equation}
which is exactly \eqref{E:HV-U}.
\end{proof}

\begin{prop}\label{P:SympL}
The symplectic leaf passing through a particular CMV matrix contains
only CMV matrices and all have the same determinant.
\end{prop}

\begin{proof}
Proving that the determinant is a Casimir is easy, so let us start
there.  Given~$\theta$, let $\phi_\theta(B)=\Im[
e^{i\theta}\log\det(B)]$ on an open neighbourhood of $\C$.  Then
$\nabla\phi_\theta= e^{i\theta} B^{-1}$ and so
$R(B\nabla\phi_\theta)=R(\nabla\phi_\theta B)=e^{-i\theta}\Id$.
Therefore,
\begin{align}
\{\phi_\theta ,\psi\} &= \tfrac{1}{2} \langle e^{-i\theta}\Id,
\nabla\psi B \rangle
        - \tfrac{1}{2}\langle e^{-i\theta}\Id , B\nabla\psi \rangle
= \tfrac{1}{2} \Im \tr([\nabla\psi, B])= 0
\end{align}
for any function $\psi$ and any angle $\theta$.

Let $\C$ be a unitary matrix and $\phi:\g\to\Reals$.  From
\eqref{E:HV-U} in Lemma~\ref{L:HV}, we see that there is an
anti-hermitian matrix $A$ so that under the Hamiltonian flow
generated by~$\phi$, $\dot \C = -\C A$. Thus $\C$ remains unitary.

Now let us restrict our attention to the case of $\C$ a CMV matrix.
We will show that under the $\phi$-flow, $\C$ remains in CMV shape,
which implies that $\C$ remains a CMV matrix by
Corollary~\ref{C:DUCK}.

By \eqref{E:HV2a} from Lemma~\ref{L:HV}, there are $L,L'\in\LL$ so
that
\begin{equation}
\dot \C = L' \C - \C L.
\end{equation}
Elementary calculations show that (left or right) multiplication by
a lower triangular matrix does not change entries above the upper
staircase; they remain zero.  Moreover, exposed entries are simply
multiplied by the corresponding diagonal entry in the lower
triangular matrix.  The first fact shows that $\dot\C$ vanishes
above the upper staircase of $\C$.  The second shows that the
exposed entries in the upper staircase obey an equation of the form
$\dot \C_{i,j}=\gamma(t) \C_{i,j}$ for some real-valued function
$\gamma$; therefore, they remain positive.

This reasoning can be transferred to the lower staircase by noting
that
\begin{equation}
 \partial_t \, \C^\dagger = - \C^\dagger \dot\C \C^\dagger = \C^\dagger L' - L \C^\dagger.
\end{equation}

This completes the proof that $\C$ remains in CMV shape and hence,
of the proposition.
\end{proof}

It remains for us to show that the symplectic leaf actually fills
out the set of CMV matrices with fixed determinant.  As this
manifold has dimension $2(n-1)$, one solution to this problem would
be to find $n-1$ functions that Poisson commute and have linearly
independent Hamiltonian vector fields.  While appearing round-about,
this approach is actually rather efficient for us because there is
just such a family of Hamiltonians that we wish to study anyway,
namely, those of the form $B\mapsto\Im\tr\{f(B)\}$, for some
polynomial $f$.   Our interest stems from their relevance to the
Ablowitz-Ladik hierarchy. We begin with the simplest abstract
properties of these Hamiltonians.

\begin{prop}\label{P:CentH}
Given a polynomial $f$, let $\phi:\g\to\Reals$ by
\begin{equation}\label{E:Phi}
\phi(B)=\Im\tr\{f(B)\}.
\end{equation}
\begin{SL}
\item $\nabla\phi\big|_{B}=f'(B)$.
\item Functions of the type \eqref{E:Phi} Poisson commute.
\item Under the flow generated by $\phi$,
\begin{equation}\label{E:LP}
 \dot B =  \bigl[B,\, \tfrac12 R\bigl(B f'(B)\bigr)\bigr]
    =\bigl[B,\, \pi_\LL\bigl(B f'(B)\bigr)] = - \bigl[B,\, \pi_\A\bigl(B f'(B)\bigr)\bigr].
\end{equation}
\item There is a unique factorization
\begin{equation}\label{E:LQ}
    \exp\bigl\{ t B f'(B) \bigr\}=L(t)Q^{-1}(t)
\end{equation}
with $L(t)\in\LLL(n)$ and $Q(t)\in\U(n)$.
\item The integral curve $B(t)$ with $B(0)=B$ is
\begin{equation}\label{E:FactS}
    B(t) = L^{-1}(t) B L(t) = Q^{-1}(t) B Q(t).
\end{equation}
\end{SL}
\end{prop}

\begin{proof}
Part (i) follows for monomials by cycling the trace.   This extends
to polynomials by linearity.

Part (ii) is readily deduced from part (i): if
$\psi(B)=\Im\tr\{g(B)\}$, then $[B,\nabla\phi]= [B,\nabla\psi]=0$
and so $\{\phi,\psi\}=0$ follows immediately from the definition,
\eqref{E:GDD}.

The first inequality in \eqref{E:LP} follows immediately from
\eqref{E:HV} and the fact that $[B,\nabla\phi]=0$. The second and
third equalities follow from \eqref{E:HV2a}.

After taking the adjoint of both sides, \eqref{E:LQ} becomes the QR
factorization.

We will prove part (v) by direct computation.  As \eqref{E:FactS} is
certainly true when $t=0$, it suffices to check that all three
formulae obey the same differential equation.  We will just prove
this for $B_1(t) =  L^{-1}(t) B L(t)$ as $Q^{-1}(t) B Q(t)$ can be
treated identically.  From the definition,
\begin{align*}
\dot B_1(t) =  L^{-1}(t) B \dot L(t) - L^{-1}(t) \dot L(t) L^{-1}(t)
B L(t)
    = \bigl[ B_1(t) , L^{-1}(t) \dot L(t) \bigr]
\end{align*}
so by \eqref{E:LP}, it suffices to show that for every $t$, $
\pi_\LL\bigl(B_1 f'(B_1)\bigr) = L^{-1} \dot L . $ This can be
demonstrated by differentiating \eqref{E:LQ}:
\begin{align}
B f'(B) \exp\bigl\{ t B f'(B) \bigr\} = \partial_t \exp\bigl\{ t B
f'(B) \bigr\}
    = \dot L Q^{-1}  - L Q^{-1}\dot Q Q^{-1},
\end{align}
which implies $B f'(B) L Q^{-1} = \dot L Q^{-1}  - L Q^{-1}\dot Q
Q^{-1}$ and so
\begin{align}
B_1 f'\bigl(B_1\bigr) = L^{-1} Bf'(B) L = L^{-1} \dot L(t)  -
Q^{-1}\dot Q.
\end{align}
This shows not only that $\pi_\LL\bigl(B_1 f'(B_1)\bigr) = L^{-1}
\dot L$, but also $\pi_\A\bigl(B_1 f'(B_1)\bigr) = -Q^{-1} \dot Q$.
\end{proof}

Part (iii) of this proposition gives Lax pair representations for
the flows generated by Hamiltonians of this form.  From the abstract
theory, this is to be expected: $B\mapsto\Im\tr\{f(B)\}$ are central
functions on $\g$---that is, they are constant on conjugacy
classes---cf. \cite[\S 2.12.4]{IntSysII}. It is not difficult to see
that these Lax pairs are precisely those discovered in \cite{N}.

Of course, the key to reconciling the concrete work of Nenciu with
the abstract approach described here is precisely the discovery that
CMV matrices have a Lie theoretic interpretation.

Similarly, \eqref{E:FactS} can be deduced from the abstract theory.
Because of the sparsity of the CMV matrix and the ready availability
of good implementations of the QR algorithm, this provides an easy
method for computational studies of the Ablowitz-Ladik hierarchy.
Those familiar with the work of Deift, Li, Nanda, and Tomei
\cite{DNT_SIAM83,DLNT_93} might immediately ask if we can say anything new about
diagonalization algorithms for unitary matrices; this will be discussed at the
end of Section~\ref{sAS}.

We will continue our investigation of these Hamiltonians in Section~\ref{sESM}.
By borrowing a result from there, we can prove the main result of this section:

\begin{theorem}\label{T:SympLeaf}
The manifold of CMV matrices with fixed determinant form a
symplectic leaf in the Poisson manifold $\g$.
\end{theorem}

\begin{proof}
Let us write $z_1,\ldots,z_n$ for the eigenvalues of the CMV matrix
$\C$.  By Lagrange interpolation, we can find polynomials $F_k$ so
that $F_k(0)=0$ and $F_k(z_j)=\delta_{jk}$.  It is then elementary
to construct polynomials $f_k$ with $F_k(z) = 2 z f_k'(z)$.

Now consider the following functions on $\g$:
\begin{equation}
 \phi_k(B) =  \Im \tr\bigl( f_k(B) \bigr).
\end{equation}
By Corollary~\ref{C:LinIndep}, $\phi_1,\ldots,\phi_{n-1}$ give rise
to linearly independent Hamiltonian vector fields.  Moreover, by
Proposition~\ref{P:CentH}, these Hamiltonians Poisson commute.
Therefore, the symplectic leaf passing through $\C$ must be of
dimension no less than $2n-2$.

From Proposition~\ref{P:SympL}, we can deduce that the dimension
must be exactly $2n-2$.  As the space of CMV matrices with fixed
determinant is path connected (it is homeomorphic to $\Disk^{n-1}$),
a simple chain-of-balls argument shows that the symplectic leaf must
exhaust this manifold.
\end{proof}

In the remainder of this section, we will describe dressing
transformations and determine how they act on CMV matrices.

We need to consider three groups $G=GL(n,\Cmplx)$, $D=G\times G$,
and $G_R=\LLL(n)\times \U^{\text{r}}(n)$.  By $\U^{\text{r}}(n)$ we
mean the group of $n\times n$ unitary matrices with the order of
multiplication reversed: $\omega_1\star \omega_2=\omega_2\omega_1$.
(We will use lower-case Greek letters for group elements.)

The significance of the group $G_R$ is that the corresponding Lie
algebra is $\g$ with the bracket given by \eqref{E:DefnRBracket}.
It is also worth noting that the map $(\lambda,\omega)\mapsto
\lambda\omega^{-1}$ defines a diffeomorphism from $G_R$ to $G$.  The
inverse mapping is given by $\beta\mapsto(\beta_+,\beta_-)$ where
$\beta_+\in\LLL(n)$ and $\beta_-\in\U(n)$ are the solution of the
factorization problem $\beta=\beta_+\beta_-^{-1}$.  This problem is
always uniquely soluble; indeed, it is essentially the QR
factorization of $\beta^\dagger$.
%

Both $G$ and $G_R$ can be regarded as subgroups (and submanifolds)
of $D$ via the following embeddings:
\begin{align*}
i :  G   \hookrightarrow D &\quad\text{by}\quad     \beta\mapsto(\beta,\beta),\quad\text{and} \\
i':G_R \hookrightarrow D &\quad\text{by}\quad
(\lambda,\omega)\mapsto(\lambda,\omega^{-1}).
\end{align*}
Moreover, the product map
$$
i \cdot i':G \times G_R \to D \quad\text{by}\quad
        (\beta,\lambda,\omega)\mapsto(\beta\lambda,\beta\omega^{-1})
$$
defines a diffeomorphism (as does $i'\!\cdot i$).  The inverse
mapping is easily seen to be $\pi\oplus\pi'$ where
\begin{align*}
\pi': D \to G_R  &\quad\text{by}\quad (\xi,\eta)\mapsto
    \bigl( [(\xi^{-1}\eta)_+]^{-1}, (\xi^{-1}\eta)_- \bigr), \quad\text{and} \\
\pi : D \to G      &\quad\text{by}\quad (\xi,\eta)\mapsto
    \xi(\xi^{-1}\eta)_+ = \eta (\xi^{-1}\eta)_-.
\end{align*}
These projections also permit us to identify the coset spaces
$D/G_R$ and $G\backslash D$ with $G$ and $G_R$ respectively.

The group $D$ acts by right multiplication on $D/G_R$ and hence, by
the identification just noted, on $G$.  As the product map $i \cdot
i'$ is onto, we can get a complete understanding of this group
action by studying its restriction to the subgroups $i(G)$ and
$i'(G_R)$.

The first subgroup, $i(G)$, leads to nothing new, just the action of
$G$ on itself by left multiplication; however $i'(G_R)$ is
different, it leads to dressing transformations:
$(\lambda,\omega)\in G_R$ maps $G$ to itself via
\begin{equation}\label{E:DT}
\beta \mapsto \pi(\lambda\beta,\omega^{-1}\beta) = \lambda \beta
(\beta^{-1}\lambda^{-1}\omega^{-1}\beta)_+ =
\omega^{-1}\beta(\beta^{-1}\lambda^{-1}\omega^{-1}\beta)_-.
\end{equation}

\begin{theorem}\label{T:Orbit}
The orbit of a CMV matrix under the dressing transformations
\eqref{E:DT} is precisely the set of CMV matrices with the same
determinant.
\end{theorem}

\begin{proof}
As we have already proved Theorem~\ref{T:SympLeaf}, this follows
from the general theory of Poisson Lie groups (cf.
\cite[\S2.12]{IntSysII}).  However, as it is possible to give a
quick and concrete proof that the orbit coincides with the
symplectic leaf, we do so.

When $\beta$ is unitary, the formula for the action can be
simplified considerably:
$(\beta^{-1}\lambda^{-1}\omega^{-1}\beta)_-=\beta^{-1}\omega(\beta^{-1}\lambda^{-1})_-$
and so
\begin{equation}\label{E:DT2}
(\lambda,\omega) : \beta \mapsto (\beta^{-1}\lambda^{-1})_-.
\end{equation}
In particular, $\omega$ plays no role.  This formula also shows
immediately that the orbit of a unitary matrix contains only unitary
matrices.

As the group $G_R$ is connected, we can prove the theorem by showing
that the tangent vectors to the orbit coincide with those of the
symplectic leaf.

Writing \eqref{E:DT2} with $\lambda=e^{-tL}$ for some $L\in\LL$ and
using the fact that $\beta$ is unitary shows
\begin{equation}\label{E:DT3}
\beta \mapsto (\beta^{-1}e^{tL})_- = \beta (\beta^{-1}e^{tL}\beta)_-
    = \beta - t \beta \, \pi_\A(\beta^{-1} L \beta) + O(t^2)  .
\end{equation}
Thus, the tangent space to the orbit through $\beta$ is given by
$-\beta \, \pi_\A(\beta^{-1} L \beta)$ as $L$ varies over $\LL$.
This is exactly the tangent space to the symplectic leaf as given in
\eqref{E:HV-U} of Lemma~\ref{L:HV}.
\end{proof}

\section{The Ablowitz-Ladik Bracket}\label{s5}

As discussed earlier, the $n$-tuple $(\alpha_0,\ldots,\alpha_{n-2},\alpha_{n-1}) \in \ID^{n-1}\times
S^1$ of Verblunsky coefficients gives a system of coordinates for the manifold of CMV matrices.

As shown in Section~\ref{s4}, the determinant, $\det(\C)=(-1)^{n-1} \bar\alpha_{n-1}$ is a Casimir
and so $(\alpha_0,\ldots,\alpha_{n-2})$ can used as coordinates for the symplectic manifold of CMV
matrices with fixed determinant.   It is therefore natural to write the Gelfand--Dikij bracket in these
coordinates.  By doing so, we will recover a bracket introduced earlier for the study of the Ablowitz--Lakik
system:
\begin{align}
\{f,g\}_2 &= \sum_{j=0}^{n-2} \rho_j^2 \left[\frac{\partial
f}{\partial u_j}\frac{\partial g}{\partial v_j}-
 \frac{\partial f}{\partial v_j}\frac{\partial g}{\partial u_j}\right]
\label{PB2:uv}\\
&=2i \sum_{j=0}^{n-2} \rho_j^2 \left[\frac{\partial f}{\partial
\bar\alpha_j}\frac{\partial g}{\partial \alpha_j}-
 \frac{\partial f}{\partial \alpha_j}\frac{\partial g}{\partial
 \bar\alpha_j}\right]\label{PB2:alpha}
\end{align}
where $\alpha_j=u_j+i v_j$ for all $0\leq j\leq n-2$, and, as usual,
$$
\frac{\partial}{\partial\alpha}=\frac12
\left(\frac{\partial}{\partial u}-i\frac{\partial}{\partial
v}\right) \quad\text{and}\quad \frac{\partial}{\partial\bar\alpha}
=\frac12 \left(\frac{\partial}{\partial u}+i\frac{\partial}{\partial
v}\right).
$$
For clarity, we will call this second bracket the Ablowitz-Ladik
bracket.  Note that this differs by a factor of two from that used in \cite{N}; we will
adjust results quoted from this paper accordingly.

On the set of Verblunsky coefficients $(\alpha_0,\ldots,\alpha_{n-2})\in \ID^{n-1}$ (with fixed $\alpha_{n-1}\in S^1$)
we consider the Hamiltonians $\Re(K_m)$ and $\Im(K_m)$, where
$$
K_m=\tfrac1m \tr(\C^m)
$$
for $m\geq1$.  The evolutions under the flows generated by the real
and imaginary parts of $K_m$ in the Ablowitz-Ladik bracket were found by Nenciu~\cite{N}:
\begin{align}\label{E:Nen1}
\dot\C &= \{ \Re (K_m), \C\}_2 = - [\C,\, \pi_\A(i\C^m) ], \\
\label{E:Nen2} \dot\C &= \{ \Im (K_m), \C \}_2 = - [\C,\,
\pi_\A(\C^m) ].
\end{align}
To obtain the formulae given above from those of \cite{N}, one
should account for the difference in notation. In \cite{N}, a
subscript $+$ is used to indicate
$$
(B_+)_{ij} = \begin{cases}
    B_{ij} & : i< j \\
    \tfrac12 B_{ii} & : i=j \\
    0 & : i > j \\
\end{cases}
$$
from which it follows that $B_+ - (B_+)^\dagger = \pi_\A(B) $;
indeed both sides are antisymmetric and agree on the upper triangle.

These formulae also hold true if we replace the Ablowitz-Ladik bracket by the
Gelfand-Dikij bracket~\eqref{E:GDD}, as one can immediately see from the second identity
in \eqref{E:HV2a}. This suggests that the Ablowitz-Ladik bracket might be the one we seek.

We will indeed prove that the Gelfand--Dikij bracket defined in
Proposition~\ref{P:PS} agrees, in the $\alpha$ coordinates, with that given above.

\begin{prop}
Consider the functions $\phi_j^a(\C) = \Im[i^a \C_{jj}]$ where
$1\leq j\leq n$ and $a\in\{0,1\}$. The Gelfand-Dikij brackets among
these functions are the following:
\begin{align}
\bigl\{ \phi_j^a , \phi_{j-1}^b \big\} = -\bigl\{ \phi_{j-1}^b ,
\phi_j^a \big\}
&= \Im[ i^{b+a} \C_{j-1,j}\C_{j,j-1}] \label{79} \\
&= -\rho_{j-2}^2 \Im[ i^{b+a} \alpha_{j-3} \bar\alpha_{j-1}],
\label{80}
\end{align}
for all choices of $a$, $b$ and $j$,
\begin{equation}
\bigl\{ \phi_j^0 , \phi_j^1 \big\} = -\bigl\{ \phi_j^1 , \phi_j^0
\big\} = |\C_{j+1,j}|^2 - |\C_{j,j-1}|^2 = \rho_{j-1}^2 -
\rho_{j-2}^2, \label{81}
\end{equation}
and all other pairs Poisson commute.  Following the standard
conventions, $\alpha_{-1}=-1$, $\rho_{-1}=0$, and matrix entries
with indices laying outside the bounds of the matrix are zero.
\end{prop}

\begin{proof} By definition, $\nabla\phi_j^a=i^a E_{jj}$ where $E_{jj}$ is the elementary matrix with
a $1$ in the $(j,j)$ position.  Thus,
\begin{align}
\bigl\{ \phi_j^a , \phi_k^b \big\} &=
\tfrac12 \langle R(i^a E_{jj}\C), i^b E_{kk}\C \rangle - \tfrac12 \langle R(i^a \C E_{jj}), i^b \C E_{kk} \rangle \\
&= \Im i^b \langle e_k | \C \pi_\LL(i^a E_{jj}\C) e_k \rangle -
   \Im i^b \langle e_k | \pi_\LL(i^a \C E_{jj}) \C e_k \rangle \label{LvYUK} \\
&= - \Im i^b \langle e_k | \C \pi_\A(i^a E_{jj}\C) e_k \rangle +
    \Im i^b \langle e_k | \pi_\A(i^a \C E_{jj}) \C e_k \rangle  \label{AvYUK}
\end{align}
where the last two lines follow from $R=2\pi_\LL-\Id=\Id-2\pi_\A$
and by observing that the identity terms cancel.  The proof proceeds
differently depending on the parity of $j$.

Fix $j$, even.  Then both $A=\pi_\A(i^a E_{jj}\C)$ and $B=\pi_\A(i^a
\C E_{jj})$ have only three non-zero entries:
\begin{gather*}
A_{j,j}= i \Im (i^a \C_{j,j}), \qquad A_{j,j+1}=-\overline{A_{j+1,j}}= i^a \C_{j,j+1} \\
B_{j,j}= i \Im (i^a \C_{j,j}), \qquad
B_{j-1,j}=-\overline{B_{j,j-1}}= i^a \C_{j-1,j}
\end{gather*}
Thus by \eqref{AvYUK}, the bracket can only be non-zero for
$k\in\{j-1,j,j+1\}$.  Moreover, if $k=j-1$ then only the first
summand in \eqref{AvYUK} contributes, while for $k=j+1$, only the
second summand contributes.  Multiplying things out gives \eqref{79}
and \eqref{81}. Equation \eqref{80} and the last formula in
\eqref{81} follow from the expressions for the matrix entries.

We now treat $j$ odd using \eqref{LvYUK}.  Both $L=\pi_\LL(i^a
E_{jj}\C)$ and $M=\pi_\LL(i^a \C E_{jj})$ have only four non-zero
entries:
\begin{gather*}
L_{j,j}= \Re(i^a \C_{j,j}), \quad L_{j,j-1} = i^a \C_{j,j-1}, \quad
L_{j+1,j}= \overline{i^a \C_{j,j+1}}, \quad L_{j+2,j} = \overline{i^a \C_{j,j+2}}, \\
M_{j,j}= \Re(i^a \C_{j,j}), \quad M_{j+1,j} = i^a \C_{j+1,j}, \quad
M_{j,j-1}= \overline{i^a \C_{j-1,j}}, \quad M_{j,j-2} =
\overline{i^a \C_{j-2,j}}.
\end{gather*}
Thus the first summand in \eqref{LvYUK} is only non-zero for
$k\in\{j-1,j\}$, while the second summand contributes only for
$k\in\{j,j+1\}$.  If $k=j+1$ or $j-1$ this leads quickly to the
formulae given in the proposition.  When $k=j$ it is a little more
involved:
\begin{align*}
\bigl\{ \phi_j^a , \phi_j^b \big\} &= \Im\bigl( i^{b-a} \bigl[
|\C_{j,j+1}|^2 + |\C_{j,j+2}|^2
        - |\C_{j-1,j}|^2 - |\C_{j-2,j}|^2 \bigr] \bigr) \\
&=\Im\bigl( i^{b-a} \bigl[ - |\C_{j,j-1}|^2 + |\C_{j+1,j}|^2 \bigr]
\bigr)
\end{align*}
where we used the fact that the $j$ row and column of $\C$ are unit
vectors.  Taking $a=0$ and $b=1$ gives \eqref{81}.
\end{proof}

To simplify our calculations, we extend the Poisson brackets defined above to be bilinear over $\IC$.
Note that $\{\bar\phi,\bar\psi\}$ is then the complex conjugate of $\{\phi,\psi\}$ and so the values
of $\{\phi,\psi\}$ and $\{\phi,\bar\psi\}$ suffice to determine all Poisson brackets between the real
and imaginary parts of $\phi$ and $\psi$.

In the previous proposition, we calculated the Gelfand-Dikij brackets of the real and imaginary
parts of the diagonal entries.  After a little computation, we obtain the following equivalent
information:
\begin{gather}
\{\C_{jj},\bar\C_{jj}\} = 2i(\rho_{j-1}^2-\rho_{j-2}^2), \label{BE1}\\
\{\C_{jj},\C_{j-1,j-1}\} = -2i\rho_{j-2}^2\alpha_{j-3}\bar\alpha_{j-1}, \label{BE2} \\
\{\C_{jj},\bar\C_{j-1,j-1}\} = 0, \label{BE3}
\end{gather}
and $\C_{jj}$ commutes with all other $\C_{kk}$ and $\bar\C_{kk}$.  These equations are the key
to showing that the two brackets are the same:

\begin{theorem}\label{T:Equal}
For any $0\leq k,l\leq n-2$ the Gelfand-Dikij brackets of the
Verblunsky coefficients are given by
\begin{equation}\label{E:BrEqual}
\{\alpha_k,\alpha_l\}=0\quad\text{and}\quad\{\alpha_k,\bar\alpha_l\}=-2i\delta_{kl}\rho_k^2.
\end{equation}
That is, the Gelfand-Dikij and Ablowitz-Ladik brackets agree.
\end{theorem}

\begin{proof}
We will prove \eqref{E:BrEqual} by induction on $k+l$.
As $\C_{11}=\bar\alpha_0$, equation \eqref{BE1} settles the case $k=l=0$.
Using the $k\leftrightarrow l$ antisymmetry, we divide the inductive step into three cases:

The case $l=k$: As $\C_{k+1,k+1}=-\alpha_{k-1}\bar\alpha_k$, equation \eqref{BE1} implies
$$
-2i|\alpha_{k}|^2\rho_{k-1}^2+|\alpha_{k-1}|^2\{\bar\alpha_k,\alpha_k\}
= 2i(\rho_{k}^2-\rho_{k-1}^2),
$$
which simplifies to $\{\alpha_k,\bar\alpha_k\}=-2i\rho_k^2$, because $\alpha_{k-1}$ is non-zero
on a dense set.

The case $l>k+1$:  As $\C_{k+1,k+1}=-\alpha_{k-1}\bar\alpha_k$, with the usual convention $\alpha_{-1}=-1$,
$$
0 = \{\C_{k+1,k+1},\C_{l+1,l+1}\} = \alpha_{k-1}\alpha_{l-1} \{\bar\alpha_k,\bar\alpha_l\}.
$$
Similarly $\{\bar\C_{k+1,k+1},\C_{l+1,l+1}\}=0$ implies $\{\alpha_k,\bar\alpha_l\}=0$.

The case $l=k+1$: By \eqref{BE3}, we have $\{\C_{k+1,k+1},\bar\C_{k+2,k+2}\}=0$ from which we deduce
$\{\bar\alpha_k,\bar\alpha_{k+1}\}=0$ as above.  Finally, by \eqref{BE2} with $j=k+2$,
$$
2i\rho_{k}^2\alpha_{k-1}\bar\alpha_{k+1}
=\{\alpha_{k-1}\bar\alpha_k,\alpha_k\bar\alpha_{k+1}\}
=2i\alpha_{k-1}\bar\alpha_{k+1}\rho_{k}^2 + \alpha_{k-1}\alpha_k\{\bar\alpha_k,\bar\alpha_{k+1}\}.
$$
Cancelling gives $\{\alpha_k,\alpha_{k+1}\}=0$.
\end{proof}

\begin{remark}
The evolution of the Verblunsky coefficients in the Ablowitz-Ladik bracket under
the flow generated by $2\Re(K_1)$ is the Ablowitz-Ladik evolution
(see \cite{AL1} and \cite{AL2}):
$$
\{\alpha_j,2\Re(K_1)\}_2=i\rho_j^2(\alpha_{j-1}+\alpha_{j+1}).
$$
\end{remark}

The proof of Theorem~\ref{T:Equal} was based around the computation of
the Gelfand-Dikij brackets of diagonal entries.  While it is possible
to directly compute the Gelfand-Dikij bracket of any
pair of entries of a CMV matrix, even the formulae for the Hamiltonian
vector fields are rather unpleasant.  However, in the particular case of exposed
entries, things are not too painful:

\begin{lemma}\label{L:MEFlow}
Fix $k,l$ so that $(k,l)$ is an exposed entry for CMV matrices. The
evolution of CMV matrices under the Hamiltonian
$r_{kl}(B)=\Re(B_{kl})$ is given by
\begin{equation}\label{E:MEFlow}
 i \, \dot\C = (-1)^k r_{kl}(\C) \bigl[ E_{kk} \C - \C E_{ll} \bigr].
\end{equation}
\end{lemma}

\begin{proof}
We will use the notation $E_{kl}$ for the matrix with $1$ in the
$(k,l)$ position and zeros everywhere else.

As $\nabla r_{kl} = i E_{lk}$, equation \eqref{E:HV2a} implies
\begin{equation}\label{E:MatEntFlow}
    \dot \C = \C \pi_\LL( i E_{lk} \C ) - \pi_\LL( i \C E_{lk} ) \C.
\end{equation}
We proceed differently depending on the parity of $k$.

As noted in Definition~\ref{D:Exposed}, the exposed entries with $k$ odd
all take the form $(2j-1,2j+1)$ with $j\geq 1$.  In this case, both
$E_{lk} \C$ and $\C E_{lk}$ are lower triangular and the sole entry
on the diagonal is positive (its value is $\C_{kl}$).  Thus we are
interested in the value of $\pi_\LL(B)$ where $B$ is lower
triangular with purely imaginary diagonal.  It is easy to see that
in this case, $\pi_\LL(B)=B-D$, where $D$ is the diagonal part of
$B$. Consequently, \eqref{E:MatEntFlow} reduces to
$$
    \dot \C = -i \C_{kl} \bigl[ \C E_{ll} - E_{kk} \C \bigr]
$$
in agreement with \eqref{E:MEFlow}.

For $k$ even, $i E_{lk} \C$ and $i \C E_{lk}$ are upper triangular
with purely imaginary diagonal. For matrices $B$ of this form,
$\pi_\LL(B)=B^\dagger+D$, where $D$ is again the diagonal part of
$B$. As $\C^\dagger\C=\C\C^\dagger=\Id$ for CMV matrices, only the
contribution from $D$ survives.  The fact that this now appears with
a plus sign results in the $(-1)^k$ factor in~\eqref{E:MEFlow}.
\end{proof}

\section{Ablowitz--Ladik: Evolution of the Spectral Measure.}\label{sESM}

Shortly after the discovery of a Lax pair representation for the
Toda Lattice \cite{Flaschka}, Moser gave a complete solution for the
finite system, \cite{Moser}. Specifically, he discovered the angle
variables associated to the actions of H\'enon \cite{Henon} and
Flaschka \cite{Flaschka}.  In addition, he studied the long-time
asymptotics and determined the scattering map.

In this section, we will discuss some corresponding results for the
Ablowitz-Ladik system. The remaining parts of the analogue of
Moser's solution can be found in Section~\ref{sAS}. A special case
of Corollary~\ref{C:mudot} has already appeared, \cite{MukNak}.  The
approach used there was to begin with a special case of
\eqref{E:MuEvol} and determine the induced evolution on the
Verblunsky coefficients.

(To be precise, Moser did not check that his `angles' Poisson
commute.  This follows as a special case of \cite[Theorem 1]{DLNT_CPAM86}.)

The key observation of Moser was that the spectral measure
associated to the Jacobi matrix and the vector $e_1$ has a very
simple evolution.  We will show that the same is true for the
Ablowitz-Ladik system.  To do this, we need a few lemmas; those
eager to see the result should skip ahead to
Proposition~\ref{P:MoserSoln} and its Corollaries.

\begin{lemma}\label{L:RBP}
For any matrices $B$ and $C$,
\begin{align}
R(B^\dagger) &= R(B) + B - B^\dagger \label{E:RL1st} \\
R\bigl( BP \bigr) &= BP - P(B-B^\dagger)P \label{E:RL2nd} \\
R\bigl( [C, P ] \bigr) &= CP + PC - 2 C^\dagger P - 2
P(C-C^\dagger)P \label{E:RL3rd}
\end{align}
where $P$ is the rank-one projection $| e_1 \rangle\langle e_1 |$.
\end{lemma}

\begin{proof}
The first result is very simple: $B-B^\dagger$ is anti-Hermitian and
so $R(B-B^\dagger)=-B+B^\dagger$. The result now follows from the
fact that $R$ is linear over $\Reals$.

Next, notice that $B| e_1 \rangle\langle e_1 |$ has non-zero entries
only in the first column. Therefore, the anti-hermitian part is just
$\tfrac12 P(B-B^\dagger)P$. Equation \eqref{E:RL2nd} now follows
easily from \eqref{E:RDef}.

The value of $R(PC)$ can be deduced from \eqref{E:RL2nd} by first
applying \eqref{E:RL1st} with $B=C^\dagger P$.  The result is $R(PC)
= 2 C^\dagger P -  P C + P(C-C^\dagger)P$. Combining this with
\eqref{E:RL2nd} gives~\eqref{E:RL3rd}.
\end{proof}

\begin{lemma}\label{L:Rpsi}
Given a polynomial $g$, let $P=| e_1 \rangle\langle e_1 |$ and
$\psi:\g\to\Reals$ by
\begin{equation}\label{E:LRpsi1}
\psi(B)=\Re\, \langle e_1 | g(B) e_1 \rangle = \langle ig(B), P
\rangle
\end{equation}
then $[B,\nabla\psi] = i [ g(B), P]$ and so
\begin{equation}\label{E:LRpsi2}
R\bigl([B,\nabla\psi]\bigr) = ig(B)P + iPg(B) + 2i g(B)^\dagger P -
4i \psi(B) P.
\end{equation}
\end{lemma}

\begin{proof}
As the operations involved are linear over $\Reals$, it suffices to
treat the case $g(B)=i^l B^k$.  By cycling the trace,
$$
\nabla\psi = i^{l+1} \sum_{q=0}^{k-1} B^q P B^{k-q-1}.
$$
Thus, the sum in $[B,\nabla\psi]$ telescopes to give $i^{l+1} [ B^k
, P]$. Lastly, \eqref{E:LRpsi2} follows from this and
Lemma~\ref{L:RBP}.
\end{proof}

\begin{prop}\label{P:MoserSoln}
Given a polynomials $f$ and $g$, let us define maps from $\g$ to
$\Reals$ by
$$
\phi(B)=\Im\tr\{f(B)\} \qquad\text{and}\qquad \psi(B)=\Re\,\langle
e_1 | g(B) e_1 \rangle.
$$
Taking $\phi$ as the Hamiltonian,
\begin{equation}\label{E:psidot}
\dot \psi(B) = \Re\,\bigl\langle e_1 | B
f'(B)\bigl(g(B)+g(B)^\dagger \bigr) e_1 \bigr\rangle
    - 2\psi(B) \Re\,\langle e_1 | B f'(B) e_1 \rangle
\end{equation}
In particular, if $d\mu$ is the spectral measure associated to $e_1$
and a CMV matrix $\C$, then writing $F(z)=2 \Re zf'(z)$, we have
\begin{equation}\label{E:ESM}
\partial_t \int\! G\,d\mu = \int\! F\, G \,d\mu  -  \int\! F \,d\mu
\int\! G \,d\mu
\end{equation}
for any function $G:S^1\to\Reals$.
\end{prop}

\begin{proof}
From Proposition~\ref{P:CentH} we have
\begin{align}
\dot \psi(B) &= \tfrac12 \Bigl\langle \bigl[B,\, R\bigl(B
f'(B)\bigr)\bigr], \, \nabla\psi \Bigr\rangle = - \tfrac12
\Bigl\langle R\bigl(B f'(B)\bigr),\, [B,\nabla\psi]
\Bigr\rangle \\
&= \tfrac12 \Bigl\langle B f'(B),\, R\bigl([B,\nabla\psi]\bigr)
\Bigr\rangle.
\end{align}
Now we can apply \eqref{E:LRpsi2} from Lemma~\ref{L:Rpsi} to deduce
\begin{align}
\dot \psi(B) &= \bigl\langle i B f'(B)\bigl(g(B)+g(B)^\dagger
\bigr),\, P \bigr\rangle
        - 2\psi(B) \bigl\langle i B f'(B),\, P \bigr\rangle,
\end{align}
which says the same thing as \eqref{E:psidot}.
\end{proof}

\begin{coro}\label{C:LinIndep}
Given a pair of polynomials $f_1$ and $f_2$, the Hamiltonians
$\phi_j(B)=\Im\tr\{f_j(B)\}$ give rise to the same integral curve
through the CMV matrix $\C$ if and only if $z\mapsto\Re[
zf_1'(z)-zf_2'(z)]$ is constant on the spectrum of $\C$.
\end{coro}

\begin{proof}
From Proposition~\ref{P:SympL}, we know that both integral curves
remain in the set of CMV matrices. As a CMV matrix is entirely
determined by its spectral measure (cf. Theorem~\ref{T:2b}), it
suffices to understand whether the two Hamiltonians give rise to the
same evolution for this measure.  From \eqref{E:ESM}, the answer is
clear: the integral curves are the same if and only if
\begin{equation}
    \int\! [F_1-F_2]\, G \,d\mu  =  \int\! [F_1-F_2] \,d\mu \int\! G \,d\mu
\end{equation}
for all functions $G:S^1\to\Reals$.  That is, if and only if
$F_1-F_2$ agrees with a constant $d\mu$-almost everywhere.  Lastly,
the support of $d\mu$ is precisely the spectrum of $\C$.
\end{proof}

Since the vector $e_1$ is cyclic for any $n\times n$ CMV matrix $\C$, the associated spectral
measure $\mu$ is supported at $n$ points. Conversely, if $\mu$ is
a probability measure on the circle which is supported at $n$ points, then
by the results in Section~\ref{s2} it is the spectral measure for $(\C,e_1)$, with
$\C$ the CMV matrix representing multiplication by $z$ in $L^2(d\mu)$.
It is natural to parameterize the spectral measure, $d\mu$,
in terms of the eigenvalues and the
mass $d\mu$ gives to them:
\begin{equation}
\int f(z) \,d\mu(z) = \sum_{j=1}^n f(e^{i\theta_j})\, \mu_j.
\end{equation}
By the observations of the previous paragraph, $\theta_j$ and $\mu_j$ are well-defined
on the manifold of CMV matrices.

By Theorem~\ref{T:2b}, one may view
$\theta_1,\ldots,\theta_{n-1},\mu_1,\ldots,\mu_{n-1}$ as a system of
coordinates on the manifold of CMV matrices with fixed determinant.
From part (ii) of Proposition~\ref{P:CentH} we know that the
functions $\theta_j$ Poisson commute. In particular, they do not
change under Hamiltonians of the form $\phi(B)=\Im\tr\{f(B)\}$.  The
evolution of $\mu_j$ is easily determined from the proposition
above.

\begin{coro}\label{C:mudot}
Under the flow generated by $\phi(B)=\Im\tr\{f(B)\}$,
\begin{equation}\label{E:CCC}
\partial_t \,\log[\mu_j] = \{\phi, \log[\mu_j] \} = F(e^{i\theta_j}) - \sum_{l=1}^n F(e^{i\theta_l})\mu_l
\end{equation}
where $F(z) = 2 \Re z f'(z)$.  Consequently,
\begin{equation}\label{E:MuEvol}
    \mu_j (t)  = \frac{\exp[F(e^{i\theta_j})\,t]\,\mu_j(0)}{\sum \exp[ F(e^{i\theta_l})\,t ]\,\mu_l(0)  }
\end{equation}
and for any $j,l\in\{1,\ldots,n-1\}$,
\begin{equation}\label{E:CanCom}
    \{ \theta_l , \tfrac12 \log[\mu_j / \mu_n ] \} = \delta_{jl}.
\end{equation}
\end{coro}

\begin{proof}
Equation~\eqref{E:CCC} amounts to nothing more that \eqref{E:ESM}
when $G(z)$ is the characteristic function of the set
$\{e^{i\theta_j}\}$.  Rather than bother solving this system of ODEs
directly, it is simple enough to check that \eqref{E:MuEvol} is
indeed the solution.

Consider \eqref{E:CCC} with $f(z)=z^k$ and $f(z)=iz^k$ for
$k=1,2,\ldots,n$.  Combining pairs gives
$$
\sum ik e^{ik\theta_l} \{\theta_l, \log[\mu_j] \}
    = 2ik \Bigl[ e^{ik\theta_j} - \sum e^{ik\theta_p} \mu_p \Bigr]
$$
This is a system of equations for $\{\theta_l, \log[\mu_j] \}$ that
has a unique solution (the Vandermonde determinant does not vanish):
\begin{equation}\label{E:CanCom2}
    \{ \theta_l , \log[\mu_j] \} = 2\delta_{lj} - 2\mu_l.
\end{equation}
Subtraction now gives \eqref{E:CanCom}.
\end{proof}

Equation \eqref{E:CanCom} gives an obvious candidate for a system of
canonically conjugate (or Darboux) coordinates:
$$
  \theta_1,\ldots,\theta_{n-1}, \tfrac12 \log[\mu_1/\mu_n],\ldots, \tfrac12 \log[\mu_{n-1}/\mu_n].
$$
Moreover, the natural analogues for the Toda lattice are canonically conjugate, \cite[Theorem 1]{DLNT_CPAM86}.
In Section~\ref{s8}, we will show that the variables $\log[\mu_k/\mu_n]$ do not commute and use the values
of their brackets to find a system of Darboux coordinates.

\section{Ablowitz--Ladik: Asymptotics and Scattering.}\label{sAS}

In this section, we study the asymptotics of solutions to the
equations in the Ablowitz-Ladik hierarchy. These asymptotics will be
expressed both in terms of the spectral measure $d\mu$ and the
Verblunsky coefficients. In the former case, the answer is an
immediate corollary of the results in the previous section.

Let us fix a Hamiltonian from the Ablowitz-Ladik hierarchy:
\begin{equation}\label{E:H6}
\phi(\C) = \Im \tr \bigl(f(\C)\bigr),
\end{equation}
where $f$ is a polynomial.  As in the previous section, we will use the shorthand $F(z)=2\Re[ z f'(z)]$.

We will write the spectral measure as
\begin{equation}\label{E:mu6}
\int f(z) \, d\mu(z) = \sum_{k=1}^n f(z_k) \mu_k.
\end{equation}
Naturally, $\C$, the spectral measure, and the Verblunsky coefficients evolve under the flow
generated by $\phi$; we use $\C(t)$, $\mu_k(t)$, and $\alpha_j(t)$ to denote these
quantities at time $t$.  By Proposition~\ref{P:CentH} part (ii), the
eigenvalues do not change under this flow; hence there is no need
consider $z_k(t)$.  When we omit the time dependence, we refer to the initial data.

In defining $\mu_k$ and $z_k$, we can choose any ordering we please;
however, there is a particular condition on this choice that
simplifies the formulae below.  Namely, we require that
\begin{equation}
 \lambda_1 \geq \lambda_2 \geq \cdots \geq \lambda_n,
\end{equation}
where we use the shorthand $\lambda_k = F(z_k)=2\Re[ z_k f'(z_k)]$.
Of course generically, the ordering will be strict. Note that by the continuity of $F$
this labeling is well-defined on an open set which is invariant under the flow of $\phi$.

\begin{prop}\label{P:muA}
Under the flow generated by the Hamiltonian $\phi$, the masses have the following asymptotics
\begin{equation}\label{AsymptMu}
\log[\mu_k(t)] = -(\lambda_1 - \lambda_k)t + \log\biggl[\frac{\mu_k}{\mu_1+\cdots+\mu_\nu}\biggr] +
O(e^{-at})
\end{equation}
as $t\rightarrow\infty$.  Here $\nu$ is defined by $\lambda_1=\cdots=\lambda_\nu>\lambda_{\nu+1}$ and
$a=\lambda_1-\lambda_{\nu+1}>0$.  In particular if $k>\nu$, then $\mu_k\to0$ exponentially fast.
\end{prop}

\begin{proof}
Equation \eqref{AsymptMu} is an immediate consequence of \eqref{E:MuEvol}, which says
$$
\log[\mu_k(t)]= \lambda_k t + \log[\mu_k]- \log\left[\sum_{l=1}^n e^{\lambda_l t}\mu_l \right].
$$
The definition of $\nu$ implies
\begin{align*}
\sum_{m=1}^n \mu_m e^{\lambda_m t} = e^{\lambda_1 t} (\mu_1 + \cdots + \mu_\nu)[1+O(e^{-at})].
\end{align*}
Combining these two formulae completes the proof.
\end{proof}

We now turn to studying the asymptotics of the Verblunsky
coefficients. In order to simplify the formulae that follow, we will
use multi-index notation. All our multi-indices will be ordered.
Given a multi-index $I=(i_1<\cdots<i_l)$ of length $l$, we write
\begin{align}
\mu_I = \prod_{k=1}^l \mu_{i_k},  \qquad z_I = \prod_{k=1}^l z_{i_k},
\end{align}
and we abbreviate Vandermonde determinants as follows:
\begin{align}\label{E:VanderDefn}
\Delta(\zeta_1,\ldots,\zeta_m) &= \det(\zeta_l^{k-1}) = \prod_{1\leq j < k \leq m} [\zeta_k - \zeta_j ] \\
\Delta_I &= \Delta(z_{i_1},\ldots,z_{i_l}).
\end{align}
We follow the natural convention that the Vandermonde determinant of a single number is $1$.

The key to converting asymptotic information on the measure into
asymptotics for the Verblunsky coefficients is the following lemma.

\begin{lemma}
Let $d\mu=\sum_{k=1}^n \mu_k \delta_{z_k}$ be a discrete measure on
the unit circle and let $\hat\mu(p)=\int z^p d\mu(z)$.  Then for
each $1\leq m \leq n$,
\begin{gather}\label{E:Toe1}
\det\bigl[\hat\mu(k-l)\bigr] = \sum_{I} |\Delta_I|^2 \mu_I
    = \prod_{j=0}^{m-2} \rho_j^{2(m-j-1)}, \\
\label{E:Toe2} \det\bigl[\hat\mu(k-l-1)\bigr] = \sum_I |\Delta_I|^2
\mu_I \bar z_I
    = (-1)^{m-1} \alpha_{m-1} \prod_{j=0}^{m-2} \rho_j^{2(m-j-1)},
\end{gather}
where the determinants are taken over $1\leq k,l \leq m$ and both
sums are over ordered multi-indices $I=(1\leq i_1<i_2<\cdots<i_m
\leq n)$.
\end{lemma}

\begin{proof} The proof is an elementary application of a special case of the Cauchy--Binet Formula
(see \cite[p. 9]{Gantmacher} or \cite[Lemma 36.2]{vLW}):  Suppose
$m\leq n$.  Given an $m \times n$ matrix $A$ and a $n\times n$
diagonal matrix $D$,
\begin{equation}\label{E:CBF}
\det(A D A^\dagger) = \sum_{i_1<\cdots<i_m}  \biggl\{ \bigl|
\det(A_{k,i_l}) \bigr|^2
        \prod_{r=1}^{m} D_{i_r,i_r} \biggr\}.
\end{equation}

To prove the first equalities in \eqref{E:Toe1} and \eqref{E:Toe2},
we choose
$$
A = \begin{bmatrix}
1 & 1 & \cdots & 1 \\
z_1 & z_2 & \cdots & z_n  \\
\vdots & \vdots & & \vdots \\
z_1^{m-1} & z_2^{m-1} & \cdots & z_n^{m-1}
\end{bmatrix}
$$
with $D=\diag[\mu_1,\ldots,\mu_n]$ and $D=\diag[\mu_1\bar
z_1,\ldots,\mu_n\bar z_n]$, respectively. The right-hand side of
\eqref{E:CBF} reduces to the form given above because of
\eqref{E:VanderDefn}. Evaluating the matrix products $ADA^\dagger$
leads to the Toeplitz matrices given on the left-hand sides of
\eqref{E:Toe1} and \eqref{E:Toe2}.

To prove the second equality in \eqref{E:Toe1}, we apply row
operations to $A$:
$$
E A = \begin{bmatrix}
\Phi_0(z_1)     & \cdots & \Phi_0(z_n) \\
\vdots &  & \vdots \\
\Phi_{m-1}(z_1) & \cdots & \Phi_{m-1}(z_n)
\end{bmatrix}
$$
for some lower-triangluar matrix $E$ of determinant one.  The
orthogonality of $\Phi_j$ with respect to $d\mu$ shows that
$EADA^\dagger E^\dagger$ is a diagonal matrix with entries
$\|\Phi_j\|^2$.  The product of the squared norms of $\Phi_j$
reduces to \eqref{E:Toe1} via \eqref{E:PhiNorm}.

The second part of \eqref{E:Toe2} follows from \eqref{E:Toe1} by the
Heine Formula, \cite[Theorem~1.5.11]{Simon1}, which says $(-1)^{m-1}
\alpha_{m-1}$ is the ratio of the left-hand sides of \eqref{E:Toe1}
and~\eqref{E:Toe2}.
\end{proof}

\begin{theorem}\label{T:AlphaAss}
Fix $1\leq k \leq n-1$.  Let $B_k=\{ l : \lambda_l=\lambda_k\}$, let
$s(k)$ be the largest index not contained in $B_k$.  Then, as
$t\to\infty$,
\begin{equation}\label{E:AsympA}
\alpha_{k-1}(t) \to (-1)^{k-1} \bar z_J \frac{\sum |\Delta_{J\cup
I}|^2 \mu_I \bar z_I}
    {\sum |\Delta_{J\cup I}|^2 \mu_I},
\end{equation}
where $J=(1<2<\cdots<s(k))$ and both sums are over all ordered
multi-indices $I\subseteq B_k$ of length $k-s(k)$. In particular,
$|\alpha_{k-1}|\to 1$ if and only if $\lambda_k>\lambda_{k+1}$.

If all $\lambda_j$ are distinct,
\begin{equation}\label{E:AsympB}
\alpha_{k-1}(t) = (-1)^{k-1} \bar z_1\cdots \bar z_k \left[ 1 +
    \xi_{k-1} e^{-(\lambda_{k}-\lambda_{k+1})t} + O( e^{-\gamma t} )  \right]
\end{equation}
where
\begin{equation}
\xi_{k-1} = (z_k \bar z_{k+1} - 1) \frac{\mu_{k+1}}{\mu_{k}}
    \prod_{l=1}^{k-1}\left|\frac{z_{k+1}-z_l}{z_{k}-z_l}\right|^2
\end{equation}
and $\gamma>(\lambda_{k}-\lambda_{k+1})>0$.
\end{theorem}

\begin{proof}
By taking the ratio of \eqref{E:Toe1} and \eqref{E:Toe2}, we obtain
\begin{gather}\label{E:Heine}
\alpha_{k-1}(t) = (-1)^{k-1} \frac{ \sum |\Delta_L|^2 \mu_L(t) \bar
z_L }{\sum |\Delta_L|^2 \mu_L(t)}
\end{gather}
where the sums are over all ordered multi-indices, $L$, of length
$k$.

The key to proving \eqref{E:AsympA} is determining the
asymptotically dominant terms in these sums. By
Proposition~\ref{P:muA}, if $\lambda_p > \lambda_q$, then the ratio
$\mu_q(t)/\mu_p(t)$ converges to zero exponentially fast.
Therefore, the dominant terms arise from those multi-indices for
which $\lambda_{l_1}+\cdots+\lambda_{l_k}$ is maximal.  As the
$\lambda$s are decreasing, the maximal value is
$\lambda_1+\cdots+\lambda_k$ and the set of multi-indices that
achieve this value is exactly the collection of $J\cup I$ stated in
the theorem.  This shows that
\begin{align}\label{E:AAA}
\alpha_{k-1}(t) \approx (-1)^{k-1} \frac{\sum |\Delta_{J\cup I}|^2
\mu_{J\cup I}(t) \bar z_{J\cup I} }
    {\sum |\Delta_{J\cup I}|^2 \mu_{J\cup I}(t) },
\end{align}
where $A\approx B$ means that the ratio $A/B\to 1$ exponentially
fast.  (We will maintain this convention throughout of the proof.)
Note that $z_J$ can be factored out, leaving $z_I$.

By Proposition~\ref{P:muA},
\begin{align*}
\mu_{J\cup I}(t) \approx \frac{\mu_J \mu_I}{(\mu_1+\cdots+\mu_\nu)^k}
    \exp\bigl[(\lambda_1+\cdots+\lambda_k)t-k\lambda_1 t\bigr]
\end{align*}
for the multi-indices $I$ under consideration.  Substituting this
into \eqref{E:AAA} and cancelling common factors in the numerator
and denominator gives \eqref{E:AsympA}.

From \eqref{E:AsympA} we see that the limiting value of $(-1)^{k-1}
z_J \alpha_{k-1}(t)$ is a convex combination of points on the unit
circle.  This sum contains exactly one term if and only if
$\lambda_k>\lambda_{k+1}$.  Therefore under this condition,
$|\alpha_{k-1}(t)|\to 1$.  If the sum contains more than one term,
then it contains two multi-indices that differ at only one index.
Thus the limiting value of $\alpha_{k-1}(t)$ is a non-trivial convex
combination of points on the unit circle and so lies strictly inside
the unit disk.

We will now consider the case of distinct $\lambda_j$.  This amounts
to studying the asymptotics of $\mu_L(t)$ for different choices of
ordered multi-indices, $L=(l_1<\cdots<l_k)$.  By
Proposition~\ref{P:muA},
$$
\log[ \mu_L(t) ] \approx  - t \sum_{j=1}^k (\lambda_1 -
\lambda_{l_j}) + \log[\mu_L] - k \log[\mu_1]
$$
and so the question reduces to finding the largest two values of
$\lambda_{l_1}+\cdots+\lambda_{l_k}$. Because
$\lambda_1>\lambda_2>\cdots>\lambda_n$, the largest value is
$\lambda_1+\cdots+\lambda_k$; the second largest is
$\lambda_1+\cdots+\lambda_{k-1}+\lambda_{k+1}$.  All other choices
of $L$ lead to strictly smaller values for the sum. We write
$\tilde\gamma$ for the difference between the largest and third
largest values of this sum.

Returning to \eqref{E:AAA} and factoring out the dominant term in
both the numerator and the denominator gives
\begin{align*}
(-1)^{k-1}z_1\cdots z_k\alpha_{k-1} &= \frac{1+ z_k\bar z_{k+1} P_k
e^{-(\lambda_k-\lambda_{k+1})t}+O(e^{-\tilde\gamma t})}
    {1+ P_k e^{-(\lambda_k-\lambda_{k+1})t}+O(e^{-\tilde\gamma t})}\\
&=1+ (z_k\bar z_{k+1}-1)P_k
e^{-(\lambda_k-\lambda_{k+1})t}+O(e^{-\gamma t}),
\end{align*}
where $\gamma=\min\{\tilde\gamma, 2(\lambda_{k}-\lambda_{k+1})\}$ and
$$
P_k=\frac{\mu_{k+1}}{\mu_k}\left|\frac{\Delta(z_1,\ldots,z_{k-1},z_{k+1})}{\Delta(z_1,\ldots,z_k)}\right|^2.
$$
Cancelling common factors in the Vandermonde determinants leads to the formula given in the theorem.
\end{proof}

If all $\lambda_j$ are distinct, then viewed as a curve in the
disk, $\alpha_{k-1}(t)$ approaches the boundary in a fixed
non-tangential direction.  This simply amounts to the statement that
$\xi_{k-1}$ is non-zero and $\arg(\xi_{k-1})=\arg(z_k \bar z_{k+1} -
1)$ belongs to $(-\pi/2,\pi/2)$. Let us also note that the
asymptotics of $\rho_{k-1}$ are easily deduced from
\eqref{E:AsympB}:
\begin{align}
  \rho_{k-1}^2 (t) &= -2\Re(\xi_{k-1}) e^{-(\lambda_{k}-\lambda_{k+1})t} + O( e^{-\gamma t} )  \\
  &= | z_{k+1}-z_k |^2 \frac{\mu_{k+1}}{\mu_{k}}
    \prod_{l=1}^{k-1}\left|\frac{z_{k+1}-z_l}{z_{k}-z_l}\right|^2 e^{-(\lambda_{k}-\lambda_{k+1})t} + O( e^{-\gamma t} ).
\end{align}

This shows that the factors $\mathcal{L}(t)$ and $\mathcal{M}(t)$ of the CMV matrix $\C(t)$
diagonalize as $t\to\infty$ and hence so does $\C(t)$.  Moreover, the eigenvalues are ordered
by the corresponding value of $F(z)$. This is a well known phenomenon for the Toda Lattice.

When the $\lambda_j$ are not all distinct, $\C(t)$ converges to a direct sum of CMV matrices and
their adjoints.  Specifically, if
$$
  \lambda_{k-1} > \lambda_k=\cdots=\lambda_{k+m} > \lambda_{k+m+1},
$$
then $\alpha_{k-1},\ldots,\alpha_{k+m-2}$ do not approach the
boundary and $\C(\infty)$ has a non-trivial block of size $m$
beginning at row/column $k$.  If $k$ is odd, this will be a CMV
matrix; if $k$ is even, it will the adjoint of a CMV matrix.

While this phenomenon cannot occur for the Toda Lattice, it can
occur for Hamiltonians in the same hierarchy.  Some examples of
non-diagonalization are discussed on page~389 of \cite{DLT_JFA85}.
However, we do not know if anyone has troubled to write down the
analogue of Theorem~\ref{T:AlphaAss}.

We now turn our attention to scattering; we will only consider the
case of distinct $\lambda$s. Thus far, we have only discussed the
behaviour as $t\to +\infty$; however the behaviour for $t\to-\infty$
can easily be deduced from this by reversing the sign of the
Hamiltonian.  This results in the following changes to the
$t\to\infty$ formulae:
$$
\lambda_k \mapsto -\lambda_{n-k+1},\quad z_k\mapsto z_{n-k+1},
    \quad\text{and}\quad \mu_k\mapsto\mu_{n-k+1}
$$
In particular,
\begin{align*}
\alpha_{k-1}(t) = (-1)^{k-1} \bar z_n\cdots \bar z_{n-k+1} \left[ 1
+
    \zeta_{k-1} e^{(\lambda_{n-k}-\lambda_{n-k+1}) t} + O( e^{\gamma t} )  \right]
\end{align*}
as $t\to -\infty$ where
\begin{equation}
\zeta_{k-1} = (z_{n-k+1} \bar z_{n-k} - 1)
\frac{\mu_{n-k}}{\mu_{n-k+1}}
    \prod_{l=n-k+2}^{n} \left|\frac{z_{n-k}-z_l}{z_{n-k+1}-z_l}\right|^2
\end{equation}
and $\gamma> (\lambda_{n-k}-\lambda_{n-k+1}) >0 $.

The scattering map is the transformation that links the asymptotics of solutions
as $t\to-\infty$ to those as $t\to\infty$.  The result follows easily from the
formulae given above:

\begin{prop}
Suppose all $\lambda_j$ are distinct, then
$$
\alpha_{k-1}(+\infty) \alpha_{n-k-1}(-\infty) = (-1)^{n} \bar
z_1\cdots \bar z_n
    = (-1)^n \det(\C^\dagger) = -\alpha_{n-1},
$$
which is independent of $t$.  Furthermore,
$$
\xi_{k-1}\zeta_{n-k-1} = |z_{k+1} - z_{k}|^2
        \prod_{l=1}^{k-1}\left|\frac{z_{k+1}-z_l}{z_{k}-z_l}\right|^2
        \prod_{m=k+2}^{n} \left|\frac{z_{k}-z_m}{z_{k+1}-z_m}\right|^2.
$$
The eigenvalues $z_j$ are time independent and can be recovered from the asymptotic values of
the Verblunsky coefficients as follows:
$$
z_{j} = -\frac{\alpha_{j-2}(+\infty)}{\alpha_{j-1}(+\infty)} = -\frac{\alpha_{n-j-1}(-\infty)}{\alpha_{n-j}(-\infty)}.
$$
\end{prop}

The corresponding result for the Toda Lattice can be found in \cite[\S 4]{Moser}.

As noted earlier, the Toda flow diagonalizes Jacobi matrices.  One may ask if it suggests a good
algorithm for doing this.  Indeed it does, namely the venerable QR algorithm. See \cite{DLNT_93,Deift}.

For unitary matrices, the naive QR algorithm is useless: the QR factorization of $U$ is $U \cdot \Id$.
However, by incorporating shifting and deflation one does obtain an algorithm that works.  In view of
parts (iv) and (v) of Proposition~\ref{P:CentH}, any algorithm based on the Ablowitz--Ladik flows
will merely amount to the application of the QR method to a (possibly step dependent) function of
$U$.  We do not wish to be entirely pessimistic, particularly in light of the storage savings afforded
by CMV matrices.  One may try the following algorithm: Iterate
$$
\C \mapsto Q\C Q^\dagger \quad\text{where} \quad \C\pm\C^\dagger =LQ.
$$
The choice of sign may be fixed or vary.
Whether this has any real numerical benefit, is a matter for experiment.

\section{Canonical coordinates}\label{s8}

In this section, we extend $\theta_1,\ldots,\theta_{n-1}$ to a system of Darboux
coordinates; these Poisson commute by Proposition~\ref{P:CentH}.
Based on \eqref{E:CanCom} it is natural to imagine that
\begin{equation}
  \tfrac12\log[\mu_1/\mu_n],\ldots, \tfrac12 \log[\mu_{n-1}/\mu_n]
\end{equation}
are the missing variables as is true for Jacobi matrices (a special case of \cite[Theorem~1]{DLNT_CPAM86}).
As we will see, they are not, because they do not Poisson commute.  The remedy is well known
(see \cite[Ch. 4]{AKN}): the conjugate variables can be written in the form
$$
\varphi_l = \tfrac12\log[\mu_l/\mu_n] + f_l(\theta_1,\ldots,\theta_{n-1})
$$
where $f_l$ are determined by the system of differential equations
\begin{equation}\label{E:Poin}
\frac{\partial f_l}{\partial \theta_k} - \frac{\partial f_k}{\partial \theta_l}
        = \{ \tfrac12\log[\mu_k/\mu_n] , \tfrac12\log[\mu_l/\mu_n] \}.
\end{equation}
This system is soluble via the Poincar\'e lemma.

We begin by calculating the right hand side of \eqref{E:Poin}.
There doesn't seem to be any easy way to do this.  Because of the permutation symmetry in the indices,
it suffices to compute $\bigl\{ \log[\mu_2/\mu_1] , \log[\mu_3/\mu_1] \bigr\}$.

\begin{lemma}\label{L:PBconst}
For any enumeration of the eigenvalues $e^{i\theta_1},\ldots,e^{i\theta_n}$, the functions
$$
  \Psi_{j,k}(\C) = \bigl\{ \log[\mu_j/\mu_1] , \log[\mu_k/\mu_1] \bigr\}
$$
Poisson commute with all Hamiltonians of the form $\Im \tr f(\C)$
and so remain constant under the corresponding flows.
\end{lemma}

\begin{proof}
By the Leibnitz rule, it suffices to prove that $\Psi_{j,k}$
commutes with $\theta_l$ for each $l=1,\ldots,n$.  For this case, we
merely need to look to the Jacobi identity and \eqref{E:CanCom2}:
\begin{align*}
- \bigl\{ \bigl\{ &\log[\mu_j/\mu_1] , \log[\mu_k/\mu_1] \bigr\}, \theta_l \bigr\} \\
&= \bigl\{ \bigl\{ \log[\mu_k/\mu_1], \theta_l \bigr\},
\log[\mu_j/\mu_1] \bigr\} +
    \bigl\{ \bigl\{ \theta_l , \log[\mu_j/\mu_1] \bigr\}, \log[\mu_k/\mu_1] \bigr\} \\
&= \bigl\{ -2\delta_{lk} ,  \log[\mu_j/\mu_1] \bigr\} +
    \bigl\{ 2\delta_{lj} , \log[\mu_k/\mu_1] \bigr\} \\
&=0,
\end{align*}
which completes the proof.
\end{proof}

\begin{lemma}\label{L:Anal}
Fix distinct points $\zeta_1,\ldots,\zeta_n$ on the unit circle and
consider the collection of CMV matrices with determinant $\prod \zeta_k$ that lie within the
$\epsilon$-neighbourhood of $\diag(\zeta_1,\ldots,\zeta_n)$ as measured in the
Hilbert--Schmidt norm.  For $\epsilon$ sufficiently small, there are analytic functions $G$ and $G_j$
of the $3(n-1)$ variables $\alpha_k,\bar\alpha_k,\rho_k$, $0\leq k \leq n-2$, so that
\begin{align}
   \mu_1 &= 1 - \rho_0^2 G(\boldsymbol\alpha,\bar{\boldsymbol\alpha},\boldsymbol\rho) \label{E:An1}\\
   \log[\tfrac{\mu_2}{\mu_1}] &= \log[\rho_0^2] + \log[\tfrac{\bar\alpha_1}{(\bar\alpha_0+\bar\alpha_1\alpha_0)^2}]
    + \sum \rho_j G_j(\boldsymbol\alpha,\bar{\boldsymbol\alpha},\boldsymbol\rho), \label{E:An2}
\end{align}
where $\mu_j$ is the mass associated to the eigenvalue nearest $\zeta_j$.
\end{lemma}

\begin{proof}
Let $D$ denote the diagonal of the CMV matrix and $E$, the
remainder: $\C=D+E$.  We will also write $P_j$ for the orthogonal
projection onto the space spanned by $e_j$.  These are, of course,
the eigenprojections for $D$.

As the $\zeta_j$ are distinct, we can construct disjoint circular
contours around each of these points.  We orient these in the
anti-clockwise direction and label them $\Gamma_j$.

For $\epsilon$ sufficiently small, each contour encircles exactly
one eigenvalue of $\C$ and the corresponding eigenprojection is
given by
\begin{equation*}
\begin{aligned}
\tilde P_j  = \oint_{\Gamma_j} (z-D-E)^{-1} \tfrac{dz}{2\pi i}
= P_j &+ \oint_{\Gamma_j} (z-D)^{-1}E(z-D)^{-1} \tfrac{dz}{2\pi i}\\
    &+\oint_{\Gamma_j} (z-D)^{-1}E(z-D-E)^{-1}E(z-D)^{-1} \tfrac{dz}{2\pi i}.
\end{aligned}
\end{equation*}
The second equality follows from two applications of the resolvent
identity.

Now $\mu_j=\tr(P_1 \tilde P_j) = \langle e_1 | \tilde P_j e_1
\rangle$ and so using the fact that $D_{11}=\bar\alpha_0$, we obtain
$$
\mu_j= \delta_{1j} + \tfrac1{2\pi i}\oint_{\Gamma_j}  \langle e_1 | E (z-D-E)^{-1}E e_1 \rangle
     \frac{dz}{(z-\bar\alpha_0)^2}.
$$
The term that is linear in $E$ disappears because it has only a
double pole when $j=1$ and is analytic when $j\neq 1$.  As $Ee_1=\rho_0e_2$ and
$E^\dagger e_1=\alpha_1\rho_0e_2+\rho_0\rho_1e_3$,
\begin{equation}\label{AnalBlah}
\mu_j= \delta_{1j} + \tfrac{\rho_0^2}{2\pi i} \oint_{\Gamma_j}  \langle \alpha_1e_2+\rho_1e_3 | (z-D-E)^{-1}
    e_2 \rangle  \frac{dz}{(z-\bar\alpha_0)^2}.
\end{equation}
To see that the integral gives an analytic function of $\boldsymbol\alpha$, $\bar{\boldsymbol\alpha}$,
and $\boldsymbol\rho$, we need merely write $(z-D-E)^{-1} = \sum (z-D)^{-1} [E (z-D)^{-1}]^l$,
which is norm convergent in a neighbourhood of $\Gamma_j$ provided $\|E\|<\epsilon$ is
sufficiently small.  Taking $j=1$ immediately gives \eqref{E:An1}.

Elementary manipulations show that \eqref{E:An2} follows from
\begin{equation}\label{AnMu2}
\mu_2 = \rho_0^2\tfrac{\bar\alpha_1}{(\bar\alpha_0+\bar\alpha_1\alpha_0)^2}
    + \rho_0^2 \sum \rho_j \tilde G_j(\boldsymbol\alpha,\bar{\boldsymbol\alpha},\boldsymbol\rho),
\end{equation}
for some analytic functions $\tilde G_j$.  To prove this, we simply need to apply the resolvent formula
one more time. From \eqref{AnalBlah} and $D_{22}=-\alpha_0\bar\alpha_1$,
\begin{equation*}
\mu_2 = \frac{\rho_0^2\bar\alpha_1}{(\alpha_0\bar\alpha_1+\bar\alpha_0)^2}
+ \rho_0^2 \oint_{\Gamma_2}  \frac{ \langle \alpha_1e_2+\rho_1e_3 |
        (z-D-E)^{-1} E e_2 \rangle \ dz}{2\pi i (z+\alpha_0\bar\alpha_1)(z-\bar\alpha_0)^2}.
\end{equation*}
As $Ee_2= \rho_0\bar\alpha_1 e_1 + \rho_1\bar\alpha_2 e_3 + \rho_1\rho_2 e_4$, an additional factor of
$\rho_0$ or $\rho_1$ can be extracted from the integral to give \eqref{AnMu2}.
\end{proof}

\begin{lemma}\label{L:Anal2}
Under the hypotheses of Lemma~\ref{L:Anal}, there exist analytic functions $H_j$
of the $3(n-1)$ variables $\alpha_k,\bar\alpha_k,\rho_k$, $0\leq k \leq n-2$, so that
\begin{align}
   \log[\tfrac{\mu_3}{\mu_1}] &= \log[\rho_0^2\rho_1^2]
    + \log\biggl[\frac{-\alpha_1\bar\alpha_2^2}{\Phi_2(-\alpha_1\bar\alpha_2)^2}\biggr]
    + \sum \rho_j H_j(\boldsymbol\alpha,\bar{\boldsymbol\alpha},\boldsymbol\rho),
\end{align}
where (consistent with the notation of Section~\ref{s2})
$$
\Phi_2(z)=z^2 + (\alpha_0\bar\alpha_1 - \bar\alpha_0)z - \bar\alpha_1
$$
and $\mu_j$ is the mass associated to the eigenvalue nearest $\zeta_j$.
\end{lemma}

\begin{proof}
In view of \eqref{E:An1}, it suffices to show that
\begin{equation}\label{E:MuAn1}
\mu_3 = \rho_0^2\rho_1^2 \biggl[ \frac{-\alpha_1\bar\alpha_2^2}{\Phi_2(-\alpha_1\bar\alpha_2)^2}
    + \sum \rho_j \tilde H_j(\boldsymbol\alpha,\bar{\boldsymbol\alpha},\boldsymbol\rho) \biggr],
\end{equation}
for some analytic functions $\tilde H_j$. This result does not follow by using the resolvent formula
once more in the previous proof; instead, we will use a block decomposition.

Let us write $\C=\bigl[ \begin{smallmatrix} A & B \\ C & D+ E\end{smallmatrix} \bigr]$ where
\begin{equation*}
A = \begin{bmatrix} \bar\alpha_0 & \rho_0\bar\alpha_1 \\ \rho_0 & -\alpha_0\bar\alpha_1 \end{bmatrix},\ %
B = \begin{bmatrix} \rho_0\rho_1 & 0 & \cdots \\ -\alpha_0\rho_1 & 0 & \cdots \end{bmatrix},\ %
C^T=\begin{bmatrix} 0 & 0 & 0 & \cdots \\ \rho_1\bar\alpha_2 & \rho_1\rho_2 & 0 & \cdots \end{bmatrix},
\end{equation*}
and $D+E$ is the splitting of the remaining block into its diagonal and off-diagonal parts.  By the well known
formulae for block matrix inversion, the top left $2\times2$ block of $(z-\C)^{-1}$ is given by
$$
(z-A)^{-1} + (z-A)^{-1}B[z-D-E-C(z-A)^{-1}B]^{-1}C(z-A)^{-1},
$$
which is relevant because as in the previous lemma,
$$
\mu_3 = \tfrac{1}{2\pi i} \oint_{\Gamma_3}  \langle e_1 | (z-\C)^{-1} e_1 \rangle  \, dz.
$$
For $\epsilon$ sufficiently small, $A$ has no eigenvalues inside $\Gamma_3$ and so the first summand
above does not contribute.

To continue, we calculate
$$
(z-A)^{-1}=\frac{1}{\Phi_2(z)}
\begin{bmatrix}z+\alpha_0\bar\alpha_1 & \rho_0\bar\alpha_1 \\
                 \rho_0 & z-\bar\alpha_0\end{bmatrix},
$$
where $\det (z-A)=\Phi_2(z)=z^2 + (\alpha_0\bar\alpha_1 - \bar\alpha_0)z - \bar\alpha_1$.
This immediately implies that
$$
C(z-A)^{-1}e_1=\frac{\rho_0\rho_1}{\Phi_2(z)} \ [\bar\alpha_2 e_1 + \rho_2 e_2]
$$
and
$$
e_1^T (z-A)^{-1}B=\frac{z\rho_0\rho_1}{\Phi_2(z)} \ e_1^T.
$$
Putting all these formulae together gives
\begin{equation*}
\mu_3 = \rho_0^2\rho_1^2 \oint_{\Gamma_3}  \langle e_1 |
            [z-D-E-C(z-A)^{-1}B]^{-1} (\bar\alpha_2 e_1 + \rho_2 e_2) \rangle  \, \tfrac{z\,dz}{2\pi i\Phi_2(z)^2}.
\end{equation*}
From here we proceed as usual: let $\tilde E=E+C(z-A)^{-1}B$ and apply the resolvent identity
$(z-D-\tilde E)^{-1}=(z-D)^{-1}+(z-D)^{-1}\tilde E(z-D-\tilde E)^{-1}$ in the equation above.
As $D_{11}=\C_{33}=-\alpha_1\bar\alpha_2$, the first term gives rise to
$$
\rho_0^2\rho_1^2\bar\alpha_2 \oint_{\Gamma_3}  \langle e_1 |
            (z-D)^{-1}  e_1 \rangle  \, \frac{z\,dz}{2\pi i\Phi_2(z)^2}
=\rho_0^2\rho_1^2\,\frac{-\alpha_1\bar\alpha_2^2}{\Phi_2(-\alpha_1\bar\alpha_2)^2}\,,
$$
which is exactly the dominant term in \eqref{E:MuAn1}.  The remaining term can be rewritten
$$
\rho_0^2\rho_1^2 \oint_{\Gamma_3}  \langle e_1 | \tilde E (z-D-\tilde E)^{-1}
            (\bar\alpha_2 e_1 + \rho_2 e_2) \rangle  \, \frac{z\,dz}{2\pi i(z+\alpha_1\bar\alpha_2)\Phi_2(z)^2}
$$
which is of the form $\rho_0^2\rho_1^2\sum \rho_j
    \tilde H_j(\boldsymbol\alpha,\bar{\boldsymbol\alpha},\boldsymbol\rho)$, because every entry of $\tilde E$
contains at least one factor $\rho_j$.
\end{proof}

\begin{prop}\label{P:OhDear}
For any labelling of the eigenvalues,
\begin{equation}\label{E:OhDear}
\bigl\{ \log[\mu_2/\mu_1] , \log[\mu_3/\mu_1] \bigr\} = 2\cot\bigl(\tfrac{\theta_1-\theta_2}2\bigr)
+ 2\cot\bigl(\tfrac{\theta_2-\theta_3}2\bigr)
+ 2\cot\bigl(\tfrac{\theta_3-\theta_1}2\bigr)
\end{equation}
in the Gelfand--Dikij (or Ablowitz--Ladik) bracket.
\end{prop}

\begin{proof}
We can choose a polynomial $f$ so that $F(z_j)$ are all distinct and
$$
F(z_1) > F(z_2) > F(z_3)
$$
are the three largest values.  (We remind the reader that $F(z)=2\Re z f'(z)$.)
By Lemma~\ref{L:PBconst}, $\{\tfrac12 \log[\mu_2/\mu_1] , \tfrac12 \log[\mu_3/\mu_1]\}$
is conserved under the flow generated by $\phi(\C)=\Im\tr f(\C)$. The key idea is to use
results from Section~\ref{sAS} to evaluate this bracket as $t\to\infty$.

For $t$ sufficiently large, we can apply Lemmas~\ref{L:Anal} and~\ref{L:Anal2} because $\C(t)$ converges
to a diagonal matrix under the $\phi$-flow.  Thus we need only show that
$$
\Bigl\{
\log[\rho_0^2] + \log\Bigl[\frac{\bar\alpha_1}{(\bar\alpha_0+\bar\alpha_1\alpha_0)^2}\Bigr]
    + \sum \rho_j G_j
,\ \log[\rho_0^2\rho_1^2]
    + \log\Bigl[\frac{-\alpha_1\bar\alpha_2^2}{\Phi_2(-\alpha_1\bar\alpha_2)^2}\Bigr]
    + \sum \rho_j H_j \Bigr\}
$$
converges to the right-hand side of \eqref{E:OhDear} as $t\to\infty$.  To do this we employ linearity,
the Leibnitz rule, and the following consequences of \eqref{PB2:alpha}:
\begin{align*}
\{ \log[\rho_m], \rho_m h(\alpha_m) \} &= \{ \rho_m, h(\alpha_m) \} = -i \rho_m h'(\alpha_m)  \to 0 \\
\{ \log[\rho_m], \rho_m h(\bar\alpha_m) \} &= \{ \rho_m, h(\bar\alpha_m) \} = i \rho_m h'(\bar\alpha_m)  \to 0 \\
\{ h(\alpha_m), g(\bar\alpha_m) \} &= -2i \rho_m^2 h'(\alpha_m), g'(\bar\alpha_m) \to 0
\end{align*}
for any analytic functions $h$ and $g$.  Of the nine terms in the original bracket, this kills off
all but two:
$$
\Bigl\{ \log[\rho_0^2] ,\ \log\Bigl[\frac{-\alpha_1\bar\alpha_2^2}{\Phi_2(-\alpha_1\bar\alpha_2)^2}\Bigr] \Bigr\}
+
\Bigl\{ \log\Bigl[\frac{\bar\alpha_1}{(\bar\alpha_0+\bar\alpha_1\alpha_0)^2}\Bigr],\
\log[\rho_0^2\rho_1^2] \Bigr\}.
$$
These can be computed using \eqref{PB2:alpha}.  The result is
$$
-4i\frac{\alpha_1\bar\alpha_2(\alpha_0\bar\alpha_1+\bar\alpha_0)}{\Phi_2(-\alpha_1\bar\alpha_2)}
-2i\frac{\alpha_0\bar\alpha_1-\bar\alpha_0}{\alpha_0\bar\alpha_1+\bar\alpha_0}.
$$
From Theorem~\ref{T:AlphaAss},
\begin{alignat*}{2}
\bar\alpha_0 &\to e^{i\theta_1},   &    -\alpha_0\bar\alpha_1 &\to e^{i\theta_2}, \\
-\alpha_1\bar\alpha_2&\to e^{i\theta_3},\qquad & \Phi_2(z)&\to(z-e^{i\theta_1})(z-e^{i\theta_2})
\end{alignat*}
as $t\to\infty$.  By using these formulae,
$$
\bigl\{ \log[\mu_2/\mu_1] , \log[\mu_3/\mu_1] \bigr\}
=4i\frac{e^{i\theta_3}(e^{i\theta_1}-e^{i\theta_2})}{(e^{i\theta_3}-e^{i\theta_1})(e^{i\theta_3}-e^{i\theta_2})}
+ 2i\frac{e^{i\theta_1}+e^{i\theta_2}}{e^{i\theta_1}-e^{i\theta_2}}.
$$
The rest is trigonometry.
\end{proof}

While the proposition is stated with fixed indices $(1,2,3)$, the arbitrariness of the labelling
immediately gives the result for any triple of distinct indices.  The result is
\begin{equation}\label{99}
\begin{aligned}
\Psi_{q,r,s} :={}& \{ \log[\mu_q/\mu_s] , \log[\mu_r/\mu_s] \}  \\
={}&  2\cot\bigl(\tfrac{\theta_q-\theta_r}2\bigr)
+ 2\cot\bigl(\tfrac{\theta_r-\theta_s}2\bigr)
+ 2\cot\bigl(\tfrac{\theta_s-\theta_q}2\bigr).
\end{aligned}
\end{equation}
Note that $\Psi_{q,r,s}$ is invariant under cyclic permutations of $(q,r,s)$ and odd under transpositions.
We adopt the convention that $\Psi_{q,r,s}$ vanishes if two indices coincide.

We will now give formulae for the functions $f_l$ discussed at the beginning of this section.
Let $e^{i\eta}=\det(\C)$ and $G(x)=2\int_0^1 t \cot(xt/2)\,dt$.  Then
\begin{equation}\label{Thefs}
f_l(\theta_1,\ldots,\theta_{n-1}) = \sum_{\substack{k=1\\k\neq l}}^{n-1} \theta_k \bigl[ G(\theta_k-\theta_l)
+ G(\theta_l+\eta-\theta_n) + G(\theta_n-\eta-\theta_k)\bigr]
\end{equation}
solve the system \eqref{E:Poin}.  Given any scalar function $g$, another solution would be $f_l
+\frac{\partial g}{\partial\theta_l}$; however, we have found no choice of $g$ that leads to a simpler
formula.  It is possible to compute the integral defining $G$, but this is not an elementary function---it
involves dilogarithms.
Verifying that \eqref{Thefs} do indeed solve \eqref{E:Poin} is unenlightening; the key observations are that
$G$ is odd, $G'$ is even, and that $xG'(x)+2G(x)=2\cot(x/2)$.

While it is not immediately obvious, Proposition~\ref{P:OhDear} does permit us to determine
$\{\mu_q,\mu_r\}$.

\begin{prop}
With $\Psi$ as in \eqref{99},
\begin{equation}
  \bigl\{ \log[\mu_q] , \log[\mu_r] \bigr\} = \sum_{k=1}^n \mu_k \Psi_{q,r,k}.
\end{equation}
\end{prop}

\begin{proof}
Let us write $M_{kl}=\{\log[\mu_k],\log[\mu_l]\}$.  Regarding $\mu_1$ as a
function of the remaining masses, the Leibnitz rule gives
\begin{align*}
\Psi_{k,l,1} =M_{kl} + \sum_{p=2}^n \frac{\mu_p}{\mu_1} M_{pl}
    + \sum_{p=2}^n \frac{\mu_p}{\mu_1} M_{kp}.
\end{align*}
Multiplying through by $\mu_k$ and summing over $k=2,\ldots,n$ gives
\begin{align*}
\sum_{k=2}^n \mu_k \Psi_{k,l,1} = \sum_{k=2}^n \mu_k M_{kl}  + (1-\mu_1) \sum_{p=2}^n \frac{\mu_p}{\mu_1} M_{pl}
 + \sum_{k,p=2}^n \frac{\mu_p\mu_k}{\mu_1} M_{kp} .
\end{align*}
The last (double) sum vanishes because $M_{kp}$ is antisymmetric and so we can deduce
\begin{align*}
 \sum_{k=2}^n \mu_k \Psi_{k,l,1} = \frac{1}{\mu_1} \sum_{p=2}^n \mu_p M_{pl} = -\{\log[\mu_1],\log[\mu_l]\}.
\end{align*}
The result now follows by relabelling indices.
\end{proof}

In an earlier paper, \cite{KN}, we determined the Jacobian of the map
\begin{equation}\label{E:COV}
(\theta_1,\mu_1,\ldots,\theta_{n-1},\mu_{n-1},\theta_n) \mapsto (u_0,v_0,\ldots,u_{n-1},v_{n-2},\phi)
\end{equation}
where $\alpha_k=u_k+iv_k$ and $\alpha_{n-1}=e^{i\phi}$.
The analogous result for Jacobi matrices is due to Dumitriu and Edelman, \cite{DumE}, which served as our guide.

We asked for a simpler, more direct derivation of the Jacobian. A recent preprint of Forrester
and Rains, \cite{ForrR}, gives a very direct solution.  Percy Deift showed us a derivation of the
Jacobi matrix result using the symplectic structure naturally associated to the Toda lattice.
We would like to close with the corresponding proof for CMV matrices.  The key idea is the following:
As we can write the underlying symplectic form in either set of variables, we can view \eqref{E:COV}
as a symplectomorphism between two concrete symplectic manifolds.  In particular, it must preserve
the Liouville volume.

\begin{coro} The Jacobian of the change of variables \eqref{E:COV} is given by
$$
\det\left[\frac{\partial (u_0,v_0,\ldots,\phi)}{\partial(\theta_1,\mu_1,\ldots,\theta_n)}\right] =
  - 2^{1-n} \frac{\rho_0^2\cdots\rho_{n-2}^2}{\mu_1\cdots\mu_n}.
$$
\end{coro}

\begin{proof}
The whole proof amounts to writing the symplectic volume in each set of coordinates.  In terms of
the Verblunsky coefficients, the answer can be read off from \eqref{PB2:uv}:
\begin{equation}
\omega\wedge\cdots\wedge\omega = \left[\prod_{k=0}^{n-2} \rho_k^{-2}\right]
        du_0\wedge dv_0 \wedge\cdots\wedge du_{n-2}\wedge dv_{n-2}.
\end{equation}

To determine the volume in terms of $\theta_k$ and $\mu_k$ we begin by using the Darboux coordinates
$\theta_k,\phi_k$ constructed in this section:
\begin{align*}
\omega\wedge\cdots\wedge\omega
    &= d\theta_1\wedge d\phi_1 \wedge\cdots\wedge d\theta_{n-1}\wedge d\phi_{n-1} \\
    &=\det(J) \, d\theta_1\wedge d\mu_1 \wedge\cdots\wedge d\theta_{n-1}\wedge d\mu_{n-1}
\end{align*}
where $J$ is the $2(n-1)\times2(n-1)$ block matrix
$$
J= \begin{bmatrix} \Id & 0 \\
    \frac{\partial f_k}{\partial \theta_l} & K \end{bmatrix}
\quad\text{with}\quad
    K_{kl} = \tfrac{\partial}{\partial \mu_l} \tfrac12\log[\mu_k/\mu_n]
        = \tfrac12\mu_k^{-1} \delta_{kl} + \tfrac12\mu_n^{-1}.
$$
It is not difficult to evaluate $\det(K)$ by Gaussian elimination: first subtract the first column from all
other columns, then for each $2\leq k\leq n-1$, subtract $\mu_k/\mu_n$ times the $k$th column from the first.
The result is an upper triangular matrix whose diagonal is
$
[\frac{1}{2\mu_1}+\frac{1}{2\mu_n}+\sum_{k=2}^{n-1} \frac{\mu_k}{2\mu_1\mu_n},\,
\frac1{2\mu_2},\ldots,\frac1{2\mu_{n-1}}].
$
Recalling that $\mu_1+\cdots+\mu_n=1$, the product simplifies to give
\begin{align*}
\omega\wedge\cdots\wedge\omega &=  2^{1-n}
    \left[\prod_{k=1}^{n} \mu_k^{-1}\right] d\theta_1\wedge d\mu_1 \wedge\cdots\wedge d\theta_{n-1}\wedge d\mu_{n-1}
\end{align*}

In order to add the last coordinate, recall that $(-1)^{n-1}\bar\alpha_{n-1}=\det(\C)=\prod e^{i\theta_k}$.
Therefore $-d\phi=d\theta_1+\cdots +d\theta_n$ and so
\begin{align*}
& \left[\prod_{k=0}^{n-2} \rho_k^{-2}\right]
    du_0\wedge dv_0 \wedge\cdots\wedge du_{n-2}\wedge dv_{n-2} \wedge d\phi \\
={}& - 2^{1-n}\left[\prod_{k=1}^{n} \mu_k^{-1}\right]
    d\theta_1 \wedge d\mu_1 \wedge \cdots \wedge d\theta_{n-1} \wedge d\mu_{n-1} \wedge d\theta_n,
\end{align*}
which implies the claim by the change of variables formula of differential geometry.
\end{proof}

\begin{remark}
By placing the result at the end of such a long paper, we have perhaps given the impression that this
derivation of the Jacobian must be the most complicated of all.  However, all that we really used was
the commutativity of the eigenvalues and \eqref{E:CanCom}.  If one remains single-minded, these can be
derived very quickly from either definition of the bracket \eqref{E:GDD} or \eqref{PB2:alpha}.
\end{remark}


\end{document}